\documentclass[12pt]{amsart}
\usepackage[margin=1in,includefoot,footskip=26.9pt]{geometry}
\usepackage{amsthm}
\usepackage[T1]{fontenc}
\usepackage[utf8]{inputenc}
\usepackage[english]{babel}
\usepackage[mathscr]{eucal}
\usepackage{amsmath}
\usepackage{amssymb}
\usepackage{amsfonts}
\usepackage{mathtools}
\usepackage[sc,noBBpl]{mathpazo}
\usepackage{microtype}
\usepackage{needspace}
\usepackage{float}
\usepackage{xcolor}
\usepackage{etoolbox}
\usepackage{enumitem}
\usepackage{bigints}
\usepackage{relsize}
\usepackage{mathrsfs}
\usepackage{caption}

\numberwithin{equation}{section}

\usepackage{tikz}
\usepackage{cancel}
\usepackage{color}
\usepackage{array}
\usepackage{multirow}
\usepackage{tabularx}
\usepackage{booktabs}
\usetikzlibrary{fadings}
\usetikzlibrary{patterns}
\usetikzlibrary{shadows.blur}
\usetikzlibrary{shapes}

\theoremstyle{definition}
\newtheorem{defi}{Definition}[section]
\newtheorem{uwaga}{Remark}[section]
\theoremstyle{plain}
\newtheorem{tw}{Theorem}[section]
\newtheorem{stw}{Proposition}[section]

\makeatletter
\AtBeginEnvironment{proof}{\let\@addpunct\@gobble}
\makeatother

\newenvironment{proof.}{\begin{proof}[Proof of Theorem \ref{dodatek}.]}{\end{proof}}

\newcommand{\R}{\mathbb{R}}

\newcommand{\nat}{\mathbb{N}}

\newcommand{\hau}{\mathscr{H}}
\newcommand{\leb}{\mathscr{L}}

\usepackage[
  colorlinks=false,
  pdfborder={0 0 0.6},
  linkbordercolor={0.65 0.15 0.18},
  citebordercolor={0.10 0.45 0.25},
  urlbordercolor={0.10 0.30 0.65}
]{hyperref}

\title{An elementary proof of a tangent-cone identity
for rank-one planar Lipschitz maps}

\author{Igor Buchowiec}

\address{Faculty of Mathematics, Informatics and Mechanics, University
of Warsaw, ul. Stefana Banacha 2, 02-097 Warsaw, Poland}

\email{i.buchowiec@uw.edu.pl}

\date{November 2026}

\begin{document}
\selectlanguage{english}
\begin{abstract}
We give an elementary proof of a tangent-cone identity which might be known to specialists in geometric measure theory. It concerns Lipschitz mappings $g$ from the closed unit ball $B\subset\mathbb{R}^2$ to $\mathbb{R}^2$. Let $x$ be an interior point of differentiability of $g$ at which $Dg(x)$ has rank one. If the upper one-dimensional Hausdorff density of $g(B)$ at $g(x)$ is at most one, then the tangent cone to $g(B)$ at $g(x)$ coincides with the image of $Dg(x)$. This property is used in the proof of rigidity for the four-gradient problem. We provide a detailed, self-contained proof based on differentiability, connectedness and elementary estimates for one-dimensional Hausdorff measure, without using advanced methods from geometric measure theory.
\end{abstract}

\title[A tangent-cone identity for planar Lipschitz maps]
{An elementary proof of a tangent-cone identity for rank-one planar Lipschitz maps}




\maketitle
\selectlanguage{english}
\section{Introduction}

In this paper, we consider Lipschitz mappings $g$ from the closed unit ball $B\subset\mathbb{R}^2$ to $\mathbb{R}^2$. At an interior point of differentiability $x$ where $Dg(x)$ has rank one, we prove that the tangent cone to $g(B)$ at $g(x)$ coincides with the image of $Dg(x)$, provided that the upper one-dimensional Hausdorff density of $g(B)$ at $g(x)$ is at most one (see Definition~\ref{defgestosci}).

The identity is used in the rigidity proof of the four-gradient problem by Chlebík and Kirchheim~\cite{chlebik}. One step of the proof uses the fact that, under the density assumption, a vector is tangent to the image of $g$ at $g(x)$ if and only if it belongs to the image of $Dg(x)$.

To the best of our knowledge the proof of the result in this form is missing in the literature. The closest related published result we have found is Lemma~10.2 of Fang and Kolasiński~\cite{slawek}. Under a uniform lower density bound, it identifies the classical tangent cone to the support of a Radon measure with the approximate tangent cone of the measure. 

Our aim is to give a direct, self-contained proof using elementary methods in the planar case. This isolates the tangent-cone argument from its application to rigidity and makes explicit how the upper density bound excludes additional tangent directions.

To establish the tangent-cone identity, we prove two inclusions. The inclusion of the image of $Dg(x)$ in the tangent cone follows rather simply from differentiability. The main difficulty is the reverse inclusion. Its proof is based on an analysis of the lengths of suitable arcs in directional cones (see Section~\ref{dowod}). Using the arcwise connectedness of the image of $g$, we estimate the one-dimensional Hausdorff measure of connected parts of these arcs from below by their diameters and combine it with the upper density bound to obtain a result.

We end with an example showing that neither the upper density bound nor the tangent-cone identity follows from merely the differentiability and the rank-one assumption. 

\section{Preliminaries}
We start with several definitions.
%
\begin{defi}[Hausdorff measure]
    Let $M\subset\R^n$. For an integer $m\ge1$, the
    \emph{$m$-dimensional Hausdorff measure of $M$} is defined by
    \[
    \hau^m(M)\coloneqq\lim_{\delta\to0^+}\hau^m_\delta(M),
    \]
    where
    \[
    \hau^m_\delta(M)\coloneqq
    \inf\left\{\sum_{j=1}^{\infty}\frac{\alpha(m)}{2^m}
    (\operatorname{diam}S_j)^m:
    \begin{gathered}
      M\subseteq\bigcup_{j=1}^{\infty}S_j,\\
      S_j\subset\R^n,\quad \operatorname{diam}S_j\le\delta
    \end{gathered}\right\}.
    \]
    Here $\alpha(m)\coloneqq\leb^m(B(0,1))$, with $B(0,1)$ the
    open unit ball in $\R^m$, and the diameter bound holds for every $j$.
\end{defi}

\begin{defi}[Upper Hausdorff density]\label{defgestosci}
    For each $x\in\R^2$ the \emph{one-dimensional upper Hausdorff density of a set $M\subset\R^2$ at the point $x$} is defined as the quantity
    \[\limsup_{r\to 0^+}{\frac{1}{2r}\hau^1\left(\overline{B(x,r)}\cap M\right)}\]
    and we denote it by $\theta_1^*(M,x)$.
\end{defi}

%
\begin{defi}[Tangent vector]\label{styczny}
    Let $S\subset\R^2$ and $a\in\R^2$. A vector $v\in\R^2$ is called a \emph{tangent vector to $S$ at the point $a$}, if for each $\varepsilon>0$ there exist $x\in S$ and $r>0$ satisfying
    \[\|x-a\|<\varepsilon\text{ and }\|r(x-a)-v\|<\varepsilon.\]
    The set of all tangent vectors to $S$ at $a$ is called the \emph{tangent cone to $S$ at $a$} and is denoted by $\operatorname{Tan}(S,a)$.
\end{defi}

\begin{defi}[Directional cone]
    Let $v,a\in\R^2$ and $\varepsilon>0$. The \emph{directional cone at $a$ in direction $v$ with parameter $\varepsilon$} is
    \[
    \mathbb{E}(a,v,\varepsilon)
    \coloneqq\{x\in\R^2:\|r(x-a)-v\|<\varepsilon
    \text{ for some }r>0\}.
    \]
    If $\|v\|=1$ and $0<\varepsilon<1$, its half-angle is
    $\arcsin\varepsilon$, and its vertex $a$ is not part of the cone.
\end{defi}
We will use the following lower bound for the Hausdorff measure of a connected set.
\begin{stw}\textnormal{(Remark 2.10.12. in \cite{federer})}\label{haus}
    Let $C\subset\R^n$ be a connected set. Then the following inequality holds
    \begin{equation*}
    \textnormal{diam}\,C\leq\hau^1(C).
    \end{equation*}
\end{stw}
%
\section{Main result and outline of the proof}

\begin{tw}\label{dodatek}
    Let $B=\overline{B(0,1)}\subset\R^2$ be the closed unit ball. Assume that $g:B\to\mathbb{R}^2$ is a Lipschitz mapping. Let $x\in B(0,1)$ be a point of differentiability of $g$ such that the image of $Dg(x)$ is one-dimensional. Moreover, assume that
    \begin{equation}\label{zaloz}
    \theta_1^*(g(B),g(x))\leq 1.
    \end{equation}
    Then
    \begin{equation}\label{imdg}
    \textnormal{im}(Dg(x))=\textnormal{Tan}(g(B),g(x)).
    \end{equation}
\end{tw}
\begin{uwaga}
    Since for every $x\in B(0,1)$ we have $\textnormal{Tan}(B,x)=\R^2$,
    equality \eqref{imdg} can equivalently be written as
    \[Dg(x)(\textnormal{Tan}(B,x))=\textnormal{Tan}(g(B),g(x)).\]
\end{uwaga}

\vspace{6pt}

\subsection{Outline of the proof}

Let us briefly describe the idea of the proof.
The inclusion $\subseteq$ follows by a direct calculation using differentiability of $g$ at $x$ and does not require the density assumption.

The inclusion $\supseteq$ is proved by contradiction: we assume the existence of a tangent vector outside $\operatorname{im}Dg(x)$. The reasoning is carried out in three steps.

In the first of these, using Taylor’s formula, we show that \emph{along} the image of $Dg(x)$ we can always find points in $g(B)$ regardless of how close we get to point $g(x)$. The images of short line segments through $x$ then give a contribution of at least $1$ to the density at every sufficiently small scale. Next, there remain three situations to consider (see Figures~\ref{f1}--\ref{f3}). In each of them we show that the assumption of the existence of a tangent vector outside the image of $Dg(x)$ leads to a nonzero contribution to the upper Hausdorff density "outside the image of $Dg(x)$". This is how we handle the case presented in Figure~\ref{f3}. The first two cases can be treated more simply, since the upper density there is bounded from below by $5/4$. Thus, in each of the cases, the aforementioned upper Hausdorff density must exceed $1$, which contradicts assumption \eqref{zaloz} of the theorem. Although the argument in Step~3 also applies to the first two cases, we treat them separately to make the geometry of these configurations explicit.
\Needspace{12\baselineskip}
\section{Proof of the main result}\label{dowod}
\begin{proof.}
    We divide the reasoning into four steps.
    \begin{center}
    \textbf{Step 0. The inclusion $\subseteq$}
    \end{center}
    We begin with the simpler inclusion, $\subseteq$. The proof follows the argument in Federer \cite[3.1.21]{federer}. By translating the domain and the target, we may assume that $x=0$ and $g(0)=0$. We continue to write $B$ for the translated closed ball. Its center may have changed, but $0$ is still an interior point.

    The case $v=0$ is immediate, since $0=g(0)\in g(B)$. Let now $v\in\R^2=\operatorname{Tan}(B,0)$ be nonzero. By the definition of a tangent vector, there exist $x_j\in B$ and $r_j>0$ such that $\|x_j\|\to0$ and $\|r_jx_j-v\|\to0$ as $j\to\infty$. Since $v\ne0$, we may assume $x_j\ne0$ for every $j$.

    Using linearity of $Dg(0)$ and the equality $g(0)=0$, we can write
    \[
    \begin{aligned}
    Dg(0)(v)-r_jg(x_j)
    &=Dg(0)(v-r_jx_j)+r_j\left(Dg(0)(x_j)-g(x_j)\right)\\
    &=Dg(0)(v-r_jx_j)
      +\|r_jx_j\|\frac{Dg(0)(x_j)-g(x_j)}{\|x_j\|}.
    \end{aligned}
    \]
    Hence
    \[
    \|Dg(0)(v)-r_jg(x_j)\|
    \le \|Dg(0)(v-r_jx_j)\|
       +\|r_jx_j\|\frac{\|Dg(0)(x_j)-g(x_j)\|}{\|x_j\|}.
    \]
    The first term tends to zero by continuity of the linear map $Dg(0)$. In the second term, the sequence $\|r_jx_j\|$ is bounded because $r_jx_j\to v$, while the quotient tends to zero by differentiability of $g$ at $0$. Therefore
    \begin{equation}\label{pochodna}
    \lim_{j\to\infty}\|Dg(0)(v)-r_jg(x_j)\|=0.
    \end{equation}

    This yields $Dg(0)(v)\in\operatorname{Tan}(g(B),0)$. To verify the definition explicitly, fix $\varepsilon>0$. Since $x_j\to0$ and $g$ is continuous, there exists $j_0\in\nat$ such that $\|g(x_j)\|<\varepsilon$ for every $j\ge j_0$. By~\eqref{pochodna}, we can then choose $j_1\ge j_0$ for which $\|r_{j_1}g(x_{j_1})-Dg(0)(v)\|<\varepsilon$. Thus the point $g(x_{j_1})\in g(B)$ and the positive number $r_{j_1}$ satisfy both inequalities in Definition~\ref{styczny}.

    \needspace{5\baselineskip}
    \begin{center}
    \textbf{Step 1. Preparation using Taylor's formula}
    \end{center}

    We now prove the inclusion $\supseteq$. The remaining three steps give a proof by contradiction: a tangent vector outside the image of $Dg(0)$ will force the upper one-dimensional Hausdorff density to exceed $1$. We first introduce the notation used below.

    Keep the normalization $x=0$ and $g(0)=0$, and put $L\coloneqq\operatorname{im}Dg(0)$. Choose $u\in L$ with $\|u\|=1$, and choose $v$ such that $Dg(0)(v)=u$. Write $\lambda=\|v\|>0$. By Step~0, both $u$ and $-u$ belong to $\operatorname{Tan}(g(B),0)$.

    Suppose, towards a contradiction, that there exists $w\in\operatorname{Tan}(g(B),0)\setminus L$. Since the tangent cone is invariant under multiplication by positive scalars, we may assume that $\|w\|=1$. Choose $0<\delta_0<1/2$ so small that
    \begin{equation}\label{rozl}
    w\notin\overline{\mathbb E(0,u,\delta_0)
                         \cup\mathbb E(0,-u,\delta_0)}.
    \end{equation}
    This is possible because $w\notin L$. We can then choose $0<\varepsilon_0<1/2$ so that the closed cone $\overline{\mathbb E(0,w,\varepsilon_0)}$ is disjoint from the cones $\mathbb E(0,u,\delta_0)$ and $\mathbb E(0,-u,\delta_0)$. Fix $0<\varepsilon<\varepsilon_0$.

    From Taylor's formula at $0$, with remainder denoted by $R_0$, we obtain, for all sufficiently small $|r|$,
    \begin{equation}\label{taylor}
    g(rv)=Dg(0)(rv)+R_0(rv)=ru+R_0(rv),
    \quad \frac{\|R_0(rv)\|}{|r|\lambda}\xrightarrow{r\to 0}0.
    \end{equation}
    Choose $r_0>0$ so that $rv\in B$ whenever $|r|<r_0$. For \(0<|r|<r_0\), put
    \[
    \varepsilon_r\coloneqq
    \frac{\|R_0(rv)\|}{|r|\lambda}.
    \]
    Then \(\varepsilon_r\to0\) as \(r\to0\). Moreover,
    \begin{equation}\label{epsr}
    \|g(rv)-ru\|
    =
    |r|\lambda\varepsilon_r.
    \end{equation}
    Here and below, we use the notation $\partial B(a,0)=\{a\}$.
    Consequently, for every \(0<r<r_0\),
    \begin{equation}\label{waznee}
    \begin{split}
    g(rv)\in g(B)\cap
    \partial B(ru,r\lambda\varepsilon_r),\quad 
    g(-rv)\in g(B)\cap
    \partial B(-ru,r\lambda\varepsilon_{-r}).
    \end{split}
    \end{equation}
    Fix \(0<\delta<\delta_0\) and decrease \(r_0\), if necessary, so that
    \[
    \lambda\varepsilon_r<\min\{\delta,1/4\}
    \quad\text{for every }0<|r|<r_0.
    \]
    By the definition of the directional cones,
    \[
    g(tv)\in\mathbb E(0,u,\delta),\quad
    g(-tv)\in\mathbb E(0,-u,\delta)
    \qquad(0<t<r_0).
    \]
    These images lie in disjoint cones and are disjoint from the closed cone about $w$.

    Write $[0,p]$ for a chosen arc in $g(B)$ joining $0$ to $p$. Such arcs exist: $B$ is convex, so any two points of $g(B)$ are joined by the Lipschitz image of a segment in $B$, and a path between distinct points contains an arc. Take \(0<R<r_0\) and set
    \[
    q^-\coloneqq g(-Rv),
    \quad
    q^+\coloneqq g(Rv).
    \]
    By \eqref{waznee},
    \[q^-\in
    \partial B(-Ru,R\lambda\varepsilon_{-R}),
    \quad
    q^+\in
    \partial B(Ru,R\lambda\varepsilon_R).
    \]
    Choose arcs \([0,q^-]\) and \([0,q^+]\) in \(g(B)\).

    We distinguish three cases, according to whether these arcs have a nontrivial initial subarc contained in
    $\overline{\mathbb E(0,w,\varepsilon)}$.
    \begin{enumerate}
        \item\label{jeden}
        Both arcs have such a subarc. Choose $[0,q^{--}]\subset[0,q^-]$ and $[0,q^{++}]\subset[0,q^+]$ contained in $\overline{\mathbb E(0,w,\varepsilon)}$.

        \item\label{dwa}
        Exactly one of the two arcs has such a subarc. After interchanging the signs if necessary, choose $[0,q^{++}]\subset[0,q^+]$ contained in $\overline{\mathbb E(0,w,\varepsilon)}$.

        \item\label{trzy}
        Neither arc has a nontrivial initial subarc contained in $\overline{\mathbb E(0,w,\varepsilon)}$.
    \end{enumerate}
    These cases are mutually exclusive and exhaustive.
    Figures~\ref{f1}--\ref{f3} illustrate the three configurations and indicate the points used in the length estimates below.
    
    \vspace{-1.4cm}
\begin{figure}[H]
\centering
\tikzset{every picture/.style={line width=0.75pt}}

\begin{tikzpicture}[x=0.75pt,y=0.75pt,yscale=-1,xscale=1]

\draw    (10.75,152.13) -- (613,151.23) ;
\draw    (331,108.73) -- (331,151.37) ;
\draw [shift={(331,105.73)}, rotate = 90] [fill={rgb, 255:red, 0; green, 0; blue, 0 }  ][line width=0.08]  [draw opacity=0] (8.93,-4.29) -- (0,0) -- (8.93,4.29) -- cycle    ;
\draw   (522,151.5) .. controls (522,129.68) and (539.68,112) .. (561.5,112) .. controls (583.32,112) and (601,129.68) .. (601,151.5) .. controls (601,173.32) and (583.32,191) .. (561.5,191) .. controls (539.68,191) and (522,173.32) .. (522,151.5) -- cycle ;
\draw   (52,151.5) .. controls (52,129.68) and (69.68,112) .. (91.5,112) .. controls (113.32,112) and (131,129.68) .. (131,151.5) .. controls (131,173.32) and (113.32,191) .. (91.5,191) .. controls (69.68,191) and (52,173.32) .. (52,151.5) -- cycle ;
\draw [color={rgb, 255:red, 126; green, 211; blue, 33 }  ,draw opacity=1 ][line width=1.5]    (105,114.23) .. controls (334,2.23) and (335,6.23) .. (331,151.37) ;
\draw [shift={(331,151.37)}, rotate = 91.58] [color={rgb, 255:red, 126; green, 211; blue, 33 }  ,draw opacity=1 ][fill={rgb, 255:red, 126; green, 211; blue, 33 }  ,fill opacity=1 ][line width=1.5]      (0, 0) circle [x radius= 4.36, y radius= 4.36]   ;
\draw [shift={(105,114.23)}, rotate = 333.94] [color={rgb, 255:red, 126; green, 211; blue, 33 }  ,draw opacity=1 ][fill={rgb, 255:red, 126; green, 211; blue, 33 }  ,fill opacity=1 ][line width=1.5]      (0, 0) circle [x radius= 4.36, y radius= 4.36]   ;
\draw [color={rgb, 255:red, 126; green, 211; blue, 33 }  ,draw opacity=1 ][line width=1.5]    (332.2,150.16) .. controls (331.2,3.03) and (407.2,-43.97) .. (594.2,126.03) ;
\draw [shift={(594.2,126.03)}, rotate = 42.27] [color={rgb, 255:red, 126; green, 211; blue, 33 }  ,draw opacity=1 ][fill={rgb, 255:red, 126; green, 211; blue, 33 }  ,fill opacity=1 ][line width=1.5]      (0, 0) circle [x radius= 4.36, y radius= 4.36]   ;
\draw  [color={rgb, 255:red, 208; green, 2; blue, 27 }  ,draw opacity=1 ][line width=1.5]  (224.27,150.16) .. controls (224.27,90.56) and (272.6,42.23) .. (332.2,42.23) .. controls (391.81,42.23) and (440.13,90.56) .. (440.13,150.16) .. controls (440.13,209.77) and (391.81,258.09) .. (332.2,258.09) .. controls (272.6,258.09) and (224.27,209.77) .. (224.27,150.16) -- cycle ;
\draw  [color={rgb, 255:red, 208; green, 2; blue, 27 }  ,draw opacity=1 ] (420.67,150.16) .. controls (420.67,139.41) and (429.38,130.7) .. (440.13,130.7) .. controls (450.88,130.7) and (459.6,139.41) .. (459.6,150.16) .. controls (459.6,160.91) and (450.88,169.63) .. (440.13,169.63) .. controls (429.38,169.63) and (420.67,160.91) .. (420.67,150.16) -- cycle ;
\draw  [color={rgb, 255:red, 208; green, 2; blue, 27 }  ,draw opacity=1 ] (204.81,150.16) .. controls (204.81,139.41) and (213.52,130.7) .. (224.27,130.7) .. controls (235.02,130.7) and (243.74,139.41) .. (243.74,150.16) .. controls (243.74,160.91) and (235.02,169.63) .. (224.27,169.63) .. controls (213.52,169.63) and (204.81,160.91) .. (204.81,150.16) -- cycle ;
\draw  [color={rgb, 255:red, 208; green, 2; blue, 27 }  ,draw opacity=1 ] (369.67,150.36) .. controls (369.67,143.92) and (374.89,138.7) .. (381.33,138.7) .. controls (387.78,138.7) and (393,143.92) .. (393,150.36) .. controls (393,156.81) and (387.78,162.03) .. (381.33,162.03) .. controls (374.89,162.03) and (369.67,156.81) .. (369.67,150.36) -- cycle ;
\draw  [color={rgb, 255:red, 208; green, 2; blue, 27 }  ,draw opacity=1 ] (268.67,152.36) .. controls (268.67,145.92) and (273.89,140.7) .. (280.33,140.7) .. controls (286.78,140.7) and (292,145.92) .. (292,152.36) .. controls (292,158.81) and (286.78,164.03) .. (280.33,164.03) .. controls (273.89,164.03) and (268.67,158.81) .. (268.67,152.36) -- cycle ;
\draw [color={rgb, 255:red, 126; green, 211; blue, 33 }  ,draw opacity=1 ]   (451,167) .. controls (491,137) and (554.2,156.03) .. (594.2,126.03) ;
\draw [shift={(451,167)}, rotate = 323.13] [color={rgb, 255:red, 126; green, 211; blue, 33 }  ,draw opacity=1 ][fill={rgb, 255:red, 126; green, 211; blue, 33 }  ,fill opacity=1 ][line width=0.75]      (0, 0) circle [x radius= 3.35, y radius= 3.35]   ;
\draw [color={rgb, 255:red, 126; green, 211; blue, 33 }  ,draw opacity=1 ]   (393,150.36) .. controls (433,120.36) and (382,247.23) .. (451,167) ;
\draw [shift={(393,150.36)}, rotate = 323.13] [color={rgb, 255:red, 126; green, 211; blue, 33 }  ,draw opacity=1 ][fill={rgb, 255:red, 126; green, 211; blue, 33 }  ,fill opacity=1 ][line width=0.75]      (0, 0) circle [x radius= 3.35, y radius= 3.35]   ;
\draw [color={rgb, 255:red, 126; green, 211; blue, 33 }  ,draw opacity=1 ]   (105,114.23) .. controls (145,84.23) and (310,283.23) .. (243.74,150.16) ;
\draw [shift={(243.74,150.16)}, rotate = 243.53] [color={rgb, 255:red, 126; green, 211; blue, 33 }  ,draw opacity=1 ][fill={rgb, 255:red, 126; green, 211; blue, 33 }  ,fill opacity=1 ][line width=0.75]      (0, 0) circle [x radius= 3.35, y radius= 3.35]   ;
\draw [color={rgb, 255:red, 126; green, 211; blue, 33 }  ,draw opacity=1 ]   (243.74,150.16) -- (280.33,140.7) ;
\draw [shift={(280.33,140.7)}, rotate = 345.5] [color={rgb, 255:red, 126; green, 211; blue, 33 }  ,draw opacity=1 ][fill={rgb, 255:red, 126; green, 211; blue, 33 }  ,fill opacity=1 ][line width=0.75]      (0, 0) circle [x radius= 3.35, y radius= 3.35]   ;
\draw  [color={rgb, 255:red, 208; green, 2; blue, 27 }  ,draw opacity=1 ] (303,151.53) .. controls (303,147.39) and (306.36,144.03) .. (310.5,144.03) .. controls (314.64,144.03) and (318,147.39) .. (318,151.53) .. controls (318,155.67) and (314.64,159.03) .. (310.5,159.03) .. controls (306.36,159.03) and (303,155.67) .. (303,151.53) -- cycle ;
\draw  [color={rgb, 255:red, 208; green, 2; blue, 27 }  ,draw opacity=1 ] (344,151.53) .. controls (344,147.39) and (347.36,144.03) .. (351.5,144.03) .. controls (355.64,144.03) and (359,147.39) .. (359,151.53) .. controls (359,155.67) and (355.64,159.03) .. (351.5,159.03) .. controls (347.36,159.03) and (344,155.67) .. (344,151.53) -- cycle ;
\draw [color={rgb, 255:red, 126; green, 211; blue, 33 }  ,draw opacity=1 ]   (280.33,140.7) .. controls (320.33,110.7) and (270.5,189.03) .. (310.5,159.03) ;
\draw [shift={(310.5,159.03)}, rotate = 323.13] [color={rgb, 255:red, 126; green, 211; blue, 33 }  ,draw opacity=1 ][fill={rgb, 255:red, 126; green, 211; blue, 33 }  ,fill opacity=1 ][line width=0.75]      (0, 0) circle [x radius= 3.35, y radius= 3.35]   ;
\draw [color={rgb, 255:red, 126; green, 211; blue, 33 }  ,draw opacity=1 ]   (351.5,159.03) .. controls (301,220.23) and (353,180.36) .. (393,150.36) ;
\draw [shift={(351.5,159.03)}, rotate = 129.53] [color={rgb, 255:red, 126; green, 211; blue, 33 }  ,draw opacity=1 ][fill={rgb, 255:red, 126; green, 211; blue, 33 }  ,fill opacity=1 ][line width=0.75]      (0, 0) circle [x radius= 3.35, y radius= 3.35]   ;
\draw [color={rgb, 255:red, 0; green, 0; blue, 0 }  ,draw opacity=1 ][line width=0.75]    (561.5,151.5) ;
\draw [shift={(561.5,151.5)}, rotate = 0] [color={rgb, 255:red, 0; green, 0; blue, 0 }  ,draw opacity=1 ][fill={rgb, 255:red, 0; green, 0; blue, 0 }  ,fill opacity=1 ][line width=0.75]      (0, 0) circle [x radius= 3.35, y radius= 3.35]   ;
\draw [shift={(561.5,151.5)}, rotate = 0] [color={rgb, 255:red, 0; green, 0; blue, 0 }  ,draw opacity=1 ][line width=0.75]      (0, 0) circle [x radius= 3.35, y radius= 3.35]   ;
\draw [color={rgb, 255:red, 0; green, 0; blue, 0 }  ,draw opacity=1 ][line width=0.75]    (91.5,151.5) ;
\draw [shift={(91.5,151.5)}, rotate = 0] [color={rgb, 255:red, 0; green, 0; blue, 0 }  ,draw opacity=1 ][fill={rgb, 255:red, 0; green, 0; blue, 0 }  ,fill opacity=1 ][line width=0.75]      (0, 0) circle [x radius= 3.35, y radius= 3.35]   ;
\draw [shift={(91.5,151.5)}, rotate = 0] [color={rgb, 255:red, 0; green, 0; blue, 0 }  ,draw opacity=1 ][line width=0.75]      (0, 0) circle [x radius= 3.35, y radius= 3.35]   ;
\draw [color={rgb, 255:red, 208; green, 2; blue, 27 }  ,draw opacity=1 ][line width=0.75]    (440.13,150.16) ;
\draw [shift={(440.13,150.16)}, rotate = 0] [color={rgb, 255:red, 208; green, 2; blue, 27 }  ,draw opacity=1 ][fill={rgb, 255:red, 208; green, 2; blue, 27 }  ,fill opacity=1 ][line width=0.75]      (0, 0) circle [x radius= 3.35, y radius= 3.35]   ;
\draw [shift={(440.13,150.16)}, rotate = 0] [color={rgb, 255:red, 208; green, 2; blue, 27 }  ,draw opacity=1 ][line width=0.75]      (0, 0) circle [x radius= 3.35, y radius= 3.35]   ;
\draw [color={rgb, 255:red, 208; green, 2; blue, 27 }  ,draw opacity=1 ][line width=0.75]    (381.33,150.36) ;
\draw [shift={(381.33,150.36)}, rotate = 0] [color={rgb, 255:red, 208; green, 2; blue, 27 }  ,draw opacity=1 ][fill={rgb, 255:red, 208; green, 2; blue, 27 }  ,fill opacity=1 ][line width=0.75]      (0, 0) circle [x radius= 2.68, y radius= 2.68]   ;
\draw [shift={(381.33,150.36)}, rotate = 0] [color={rgb, 255:red, 208; green, 2; blue, 27 }  ,draw opacity=1 ][line width=0.75]      (0, 0) circle [x radius= 2.68, y radius= 2.68]   ;
\draw [color={rgb, 255:red, 208; green, 2; blue, 27 }  ,draw opacity=1 ][line width=0.75]    (280.33,152.36) ;
\draw [shift={(280.33,152.36)}, rotate = 0] [color={rgb, 255:red, 208; green, 2; blue, 27 }  ,draw opacity=1 ][fill={rgb, 255:red, 208; green, 2; blue, 27 }  ,fill opacity=1 ][line width=0.75]      (0, 0) circle [x radius= 2.68, y radius= 2.68]   ;
\draw [shift={(280.33,152.36)}, rotate = 0] [color={rgb, 255:red, 208; green, 2; blue, 27 }  ,draw opacity=1 ][line width=0.75]      (0, 0) circle [x radius= 2.68, y radius= 2.68]   ;
\draw [color={rgb, 255:red, 208; green, 2; blue, 27 }  ,draw opacity=1 ][line width=0.75]    (224.27,150.16) ;
\draw [shift={(224.27,150.16)}, rotate = 0] [color={rgb, 255:red, 208; green, 2; blue, 27 }  ,draw opacity=1 ][fill={rgb, 255:red, 208; green, 2; blue, 27 }  ,fill opacity=1 ][line width=0.75]      (0, 0) circle [x radius= 3.35, y radius= 3.35]   ;
\draw [shift={(224.27,150.16)}, rotate = 0] [color={rgb, 255:red, 208; green, 2; blue, 27 }  ,draw opacity=1 ][line width=0.75]      (0, 0) circle [x radius= 3.35, y radius= 3.35]   ;
\draw [color={rgb, 255:red, 208; green, 2; blue, 27 }  ,draw opacity=1 ][line width=0.75]    (351.5,151.53) ;
\draw [shift={(351.5,151.53)}, rotate = 0] [color={rgb, 255:red, 208; green, 2; blue, 27 }  ,draw opacity=1 ][fill={rgb, 255:red, 208; green, 2; blue, 27 }  ,fill opacity=1 ][line width=0.75]      (0, 0) circle [x radius= 2.01, y radius= 2.01]   ;
\draw [shift={(351.5,151.53)}, rotate = 0] [color={rgb, 255:red, 208; green, 2; blue, 27 }  ,draw opacity=1 ][line width=0.75]      (0, 0) circle [x radius= 2.01, y radius= 2.01]   ;
\draw [color={rgb, 255:red, 208; green, 2; blue, 27 }  ,draw opacity=1 ][line width=0.75]    (310.5,151.53) ;
\draw [shift={(310.5,151.53)}, rotate = 0] [color={rgb, 255:red, 208; green, 2; blue, 27 }  ,draw opacity=1 ][fill={rgb, 255:red, 208; green, 2; blue, 27 }  ,fill opacity=1 ][line width=0.75]      (0, 0) circle [x radius= 2.01, y radius= 2.01]   ;
\draw [shift={(310.5,151.53)}, rotate = 0] [color={rgb, 255:red, 208; green, 2; blue, 27 }  ,draw opacity=1 ][line width=0.75]      (0, 0) circle [x radius= 2.01, y radius= 2.01]   ;
\draw [color={rgb, 255:red, 74; green, 144; blue, 226 }  ,draw opacity=1 ][line width=0.75]    (355,45.36) ;
\draw [shift={(355,45.36)}, rotate = 0] [color={rgb, 255:red, 74; green, 144; blue, 226 }  ,draw opacity=1 ][fill={rgb, 255:red, 74; green, 144; blue, 226 }  ,fill opacity=1 ][line width=0.75]      (0, 0) circle [x radius= 3.69, y radius= 3.69]   ;
\draw [shift={(355,45.36)}, rotate = 0] [color={rgb, 255:red, 74; green, 144; blue, 226 }  ,draw opacity=1 ][line width=0.75]      (0, 0) circle [x radius= 3.69, y radius= 3.69]   ;
\draw [color={rgb, 255:red, 74; green, 144; blue, 226 }  ,draw opacity=1 ][line width=0.75]    (304,36.36) ;
\draw [shift={(304,36.36)}, rotate = 0] [color={rgb, 255:red, 74; green, 144; blue, 226 }  ,draw opacity=1 ][fill={rgb, 255:red, 74; green, 144; blue, 226 }  ,fill opacity=1 ][line width=0.75]      (0, 0) circle [x radius= 3.69, y radius= 3.69]   ;
\draw [shift={(304,36.36)}, rotate = 0] [color={rgb, 255:red, 74; green, 144; blue, 226 }  ,draw opacity=1 ][line width=0.75]      (0, 0) circle [x radius= 3.69, y radius= 3.69]   ;
\draw [color={rgb, 255:red, 74; green, 144; blue, 226 }  ,draw opacity=1 ][line width=0.75]    (243.74,150.16) ;
\draw [shift={(243.74,150.16)}, rotate = 0] [color={rgb, 255:red, 74; green, 144; blue, 226 }  ,draw opacity=1 ][fill={rgb, 255:red, 74; green, 144; blue, 226 }  ,fill opacity=1 ][line width=0.75]      (0, 0) circle [x radius= 3.69, y radius= 3.69]   ;
\draw [shift={(243.74,150.16)}, rotate = 0] [color={rgb, 255:red, 74; green, 144; blue, 226 }  ,draw opacity=1 ][line width=0.75]      (0, 0) circle [x radius= 3.69, y radius= 3.69]   ;
\draw [color={rgb, 255:red, 74; green, 144; blue, 226 }  ,draw opacity=1 ][line width=0.75]    (451,167) ;
\draw [shift={(451,167)}, rotate = 0] [color={rgb, 255:red, 74; green, 144; blue, 226 }  ,draw opacity=1 ][fill={rgb, 255:red, 74; green, 144; blue, 226 }  ,fill opacity=1 ][line width=0.75]      (0, 0) circle [x radius= 3.69, y radius= 3.69]   ;
\draw [shift={(451,167)}, rotate = 0] [color={rgb, 255:red, 74; green, 144; blue, 226 }  ,draw opacity=1 ][line width=0.75]      (0, 0) circle [x radius= 3.69, y radius= 3.69]   ;
\draw [color={rgb, 255:red, 74; green, 144; blue, 226 }  ,draw opacity=1 ][line width=0.75]    (393,150.36) ;
\draw [shift={(393,150.36)}, rotate = 0] [color={rgb, 255:red, 74; green, 144; blue, 226 }  ,draw opacity=1 ][fill={rgb, 255:red, 74; green, 144; blue, 226 }  ,fill opacity=1 ][line width=0.75]      (0, 0) circle [x radius= 3.69, y radius= 3.69]   ;
\draw [shift={(393,150.36)}, rotate = 0] [color={rgb, 255:red, 74; green, 144; blue, 226 }  ,draw opacity=1 ][line width=0.75]      (0, 0) circle [x radius= 3.69, y radius= 3.69]   ;
\draw [color={rgb, 255:red, 74; green, 144; blue, 226 }  ,draw opacity=1 ][line width=0.75]    (280.33,140.7) ;
\draw [shift={(280.33,140.7)}, rotate = 0] [color={rgb, 255:red, 74; green, 144; blue, 226 }  ,draw opacity=1 ][fill={rgb, 255:red, 74; green, 144; blue, 226 }  ,fill opacity=1 ][line width=0.75]      (0, 0) circle [x radius= 3.69, y radius= 3.69]   ;
\draw [shift={(280.33,140.7)}, rotate = 0] [color={rgb, 255:red, 74; green, 144; blue, 226 }  ,draw opacity=1 ][line width=0.75]      (0, 0) circle [x radius= 3.69, y radius= 3.69]   ;
\draw [color={rgb, 255:red, 74; green, 144; blue, 226 }  ,draw opacity=1 ][line width=0.75]    (351.5,159.03) ;
\draw [shift={(351.5,159.03)}, rotate = 0] [color={rgb, 255:red, 74; green, 144; blue, 226 }  ,draw opacity=1 ][fill={rgb, 255:red, 74; green, 144; blue, 226 }  ,fill opacity=1 ][line width=0.75]      (0, 0) circle [x radius= 3.69, y radius= 3.69]   ;
\draw [shift={(351.5,159.03)}, rotate = 0] [color={rgb, 255:red, 74; green, 144; blue, 226 }  ,draw opacity=1 ][line width=0.75]      (0, 0) circle [x radius= 3.69, y radius= 3.69]   ;
\draw [color={rgb, 255:red, 74; green, 144; blue, 226 }  ,draw opacity=1 ][line width=0.75]    (310.5,159.03) ;
\draw [shift={(310.5,159.03)}, rotate = 0] [color={rgb, 255:red, 74; green, 144; blue, 226 }  ,draw opacity=1 ][fill={rgb, 255:red, 74; green, 144; blue, 226 }  ,fill opacity=1 ][line width=0.75]      (0, 0) circle [x radius= 3.69, y radius= 3.69]   ;
\draw [shift={(310.5,159.03)}, rotate = 0] [color={rgb, 255:red, 74; green, 144; blue, 226 }  ,draw opacity=1 ][line width=0.75]      (0, 0) circle [x radius= 3.69, y radius= 3.69]   ;
\draw    (296,-1.9) -- (331,151.37) ;
\draw    (366,-1.9) -- (331,151.37) ;
\draw [color={rgb, 255:red, 74; green, 144; blue, 226 }  ,draw opacity=1 ][line width=0.75]    (105,114.23) ;
\draw [shift={(105,114.23)}, rotate = 0] [color={rgb, 255:red, 74; green, 144; blue, 226 }  ,draw opacity=1 ][fill={rgb, 255:red, 74; green, 144; blue, 226 }  ,fill opacity=1 ][line width=0.75]      (0, 0) circle [x radius= 3.69, y radius= 3.69]   ;
\draw [shift={(105,114.23)}, rotate = 0] [color={rgb, 255:red, 74; green, 144; blue, 226 }  ,draw opacity=1 ][line width=0.75]      (0, 0) circle [x radius= 3.69, y radius= 3.69]   ;
\draw [color={rgb, 255:red, 74; green, 144; blue, 226 }  ,draw opacity=1 ][line width=0.75]    (594.2,126.03) ;
\draw [shift={(594.2,126.03)}, rotate = 0] [color={rgb, 255:red, 74; green, 144; blue, 226 }  ,draw opacity=1 ][fill={rgb, 255:red, 74; green, 144; blue, 226 }  ,fill opacity=1 ][line width=0.75]      (0, 0) circle [x radius= 3.69, y radius= 3.69]   ;
\draw [shift={(594.2,126.03)}, rotate = 0] [color={rgb, 255:red, 74; green, 144; blue, 226 }  ,draw opacity=1 ][line width=0.75]      (0, 0) circle [x radius= 3.69, y radius= 3.69]   ;

\draw (325,159.4) node [anchor=north west][inner sep=0.75pt]    {$0$};
\draw (224,2.5) node [anchor=north west][inner sep=0.75pt]    {$\mathbb{E}(0 ,w ,\varepsilon )$};
\draw (75,160.4) node [anchor=north west][inner sep=0.75pt]    {$-Ru$};
\draw (551,159.4) node [anchor=north west][inner sep=0.75pt]    {$Ru$};
\draw (354,54.76) node [anchor=north west][inner sep=0.75pt]    {$\textcolor[rgb]{0.29,0.56,0.89}{q}\textcolor[rgb]{0.29,0.56,0.89}{^{++}}$};
\draw (306,11.4) node [anchor=north west][inner sep=0.75pt]    {$\textcolor[rgb]{0.29,0.56,0.89}{q}\textcolor[rgb]{0.29,0.56,0.89}{^{--}}$};
\draw (442.13,173.03) node [anchor=north west][inner sep=0.75pt]    {$\textcolor[rgb]{0.29,0.56,0.89}{q}\textcolor[rgb]{0.29,0.56,0.89}{_{1}^{+}}$};
\draw (383.33,153.76) node [anchor=north west][inner sep=0.75pt]    {$\textcolor[rgb]{0.29,0.56,0.89}{q}\textcolor[rgb]{0.29,0.56,0.89}{_{2}^{+}}$};
\draw (350,199.4) node [anchor=north west][inner sep=0.75pt]    {$\textcolor[rgb]{0.29,0.56,0.89}{q}\textcolor[rgb]{0.29,0.56,0.89}{_{4}^{+}}$};
\draw (236,119.4) node [anchor=north west][inner sep=0.75pt]    {$\textcolor[rgb]{0.29,0.56,0.89}{q}\textcolor[rgb]{0.29,0.56,0.89}{_{1}^{-}}$};
\draw (270,110.4) node [anchor=north west][inner sep=0.75pt]    {$\textcolor[rgb]{0.29,0.56,0.89}{q}\textcolor[rgb]{0.29,0.56,0.89}{_{2}^{-}}$};
\draw (295,163.4) node [anchor=north west][inner sep=0.75pt]    {$\textcolor[rgb]{0.29,0.56,0.89}{q}\textcolor[rgb]{0.29,0.56,0.89}{_{4}^{-}}$};
\draw (101,77.4) node [anchor=north west][inner sep=0.75pt]    {$\textcolor[rgb]{0.29,0.56,0.89}{q}\textcolor[rgb]{0.29,0.56,0.89}{^{-}}$};
\draw (388,255.4) node [anchor=north west][inner sep=0.75pt]    {$\textcolor[rgb]{0.82,0.01,0.11}{B(0,R_1)}$};
\draw (598,96.4) node [anchor=north west][inner sep=0.75pt]    {$\textcolor[rgb]{0.29,0.56,0.89}{q}\textcolor[rgb]{0.29,0.56,0.89}{^{+}}$};
\draw (528,205.4) node [anchor=north west][inner sep=0.75pt]    {$B(Ru,\lambda R\varepsilon_R)$};
\draw (28,205.4) node [anchor=north west][inner sep=0.75pt]    {$B(-Ru,\lambda R\varepsilon_{-R})$};

\end{tikzpicture}
\selectlanguage{english}
\caption{Case~\eqref{jeden}}
\label{f1}
\end{figure}

\begin{figure}[H]
\centering
\tikzset{every picture/.style={line width=0.75pt}}

\begin{tikzpicture}[x=0.75pt,y=0.75pt,yscale=-1,xscale=1]

\draw    (20,151.23) -- (613,151.23) ;
\draw    (331,108.73) -- (331,151.37) ;
\draw [shift={(331,105.73)}, rotate = 90] [fill={rgb, 255:red, 0; green, 0; blue, 0 }  ][line width=0.08]  [draw opacity=0] (8.93,-4.29) -- (0,0) -- (8.93,4.29) -- cycle    ;
\draw   (522,151.5) .. controls (522,129.68) and (539.68,112) .. (561.5,112) .. controls (583.32,112) and (601,129.68) .. (601,151.5) .. controls (601,173.32) and (583.32,191) .. (561.5,191) .. controls (539.68,191) and (522,173.32) .. (522,151.5) -- cycle ;
\draw   (52,151.5) .. controls (52,129.68) and (69.68,112) .. (91.5,112) .. controls (113.32,112) and (131,129.68) .. (131,151.5) .. controls (131,173.32) and (113.32,191) .. (91.5,191) .. controls (69.68,191) and (52,173.32) .. (52,151.5) -- cycle ;
\draw [color={rgb, 255:red, 126; green, 211; blue, 33 }  ,draw opacity=1 ][line width=1.5]    (105,114.23) .. controls (234,54.23) and (124,258.23) .. (331,151.37) ;
\draw [shift={(331,151.37)}, rotate = 332.69] [color={rgb, 255:red, 126; green, 211; blue, 33 }  ,draw opacity=1 ][fill={rgb, 255:red, 126; green, 211; blue, 33 }  ,fill opacity=1 ][line width=1.5]      (0, 0) circle [x radius= 4.36, y radius= 4.36]   ;
\draw [shift={(105,114.23)}, rotate = 335.06] [color={rgb, 255:red, 126; green, 211; blue, 33 }  ,draw opacity=1 ][fill={rgb, 255:red, 126; green, 211; blue, 33 }  ,fill opacity=1 ][line width=1.5]      (0, 0) circle [x radius= 4.36, y radius= 4.36]   ;
\draw [color={rgb, 255:red, 126; green, 211; blue, 33 }  ,draw opacity=1 ][line width=1.5]    (331,151.37) .. controls (330,4.23) and (356,2.23) .. (550,190.23) ;
\draw [shift={(550,190.23)}, rotate = 44.1] [color={rgb, 255:red, 126; green, 211; blue, 33 }  ,draw opacity=1 ][fill={rgb, 255:red, 126; green, 211; blue, 33 }  ,fill opacity=1 ][line width=1.5]      (0, 0) circle [x radius= 4.36, y radius= 4.36]   ;
\draw  [color={rgb, 255:red, 208; green, 2; blue, 27 }  ,draw opacity=1 ][line width=1.5]  (224.27,150.16) .. controls (224.27,90.56) and (272.6,42.23) .. (332.2,42.23) .. controls (391.81,42.23) and (440.13,90.56) .. (440.13,150.16) .. controls (440.13,209.77) and (391.81,258.09) .. (332.2,258.09) .. controls (272.6,258.09) and (224.27,209.77) .. (224.27,150.16) -- cycle ;
\draw  [color={rgb, 255:red, 208; green, 2; blue, 27 }  ,draw opacity=1 ] (420.67,150.16) .. controls (420.67,139.41) and (429.38,130.7) .. (440.13,130.7) .. controls (450.88,130.7) and (459.6,139.41) .. (459.6,150.16) .. controls (459.6,160.91) and (450.88,169.63) .. (440.13,169.63) .. controls (429.38,169.63) and (420.67,160.91) .. (420.67,150.16) -- cycle ;
\draw  [color={rgb, 255:red, 208; green, 2; blue, 27 }  ,draw opacity=1 ] (369.67,150.36) .. controls (369.67,143.92) and (374.89,138.7) .. (381.33,138.7) .. controls (387.78,138.7) and (393,143.92) .. (393,150.36) .. controls (393,156.81) and (387.78,162.03) .. (381.33,162.03) .. controls (374.89,162.03) and (369.67,156.81) .. (369.67,150.36) -- cycle ;
\draw [color={rgb, 255:red, 126; green, 211; blue, 33 }  ,draw opacity=1 ]   (451,167) .. controls (491,137) and (511,235.23) .. (550,190.23) ;
\draw [shift={(451,167)}, rotate = 323.13] [color={rgb, 255:red, 126; green, 211; blue, 33 }  ,draw opacity=1 ][fill={rgb, 255:red, 126; green, 211; blue, 33 }  ,fill opacity=1 ][line width=0.75]      (0, 0) circle [x radius= 3.35, y radius= 3.35]   ;
\draw [color={rgb, 255:red, 126; green, 211; blue, 33 }  ,draw opacity=1 ]   (393,150.36) .. controls (433,120.36) and (382,247.23) .. (451,167) ;
\draw [shift={(393,150.36)}, rotate = 323.13] [color={rgb, 255:red, 126; green, 211; blue, 33 }  ,draw opacity=1 ][fill={rgb, 255:red, 126; green, 211; blue, 33 }  ,fill opacity=1 ][line width=0.75]      (0, 0) circle [x radius= 3.35, y radius= 3.35]   ;
\draw [color={rgb, 255:red, 0; green, 0; blue, 0 }  ,draw opacity=1 ][line width=0.75]    (561.5,151.5) ;
\draw [shift={(561.5,151.5)}, rotate = 0] [color={rgb, 255:red, 0; green, 0; blue, 0 }  ,draw opacity=1 ][fill={rgb, 255:red, 0; green, 0; blue, 0 }  ,fill opacity=1 ][line width=0.75]      (0, 0) circle [x radius= 3.35, y radius= 3.35]   ;
\draw [shift={(561.5,151.5)}, rotate = 0] [color={rgb, 255:red, 0; green, 0; blue, 0 }  ,draw opacity=1 ][line width=0.75]      (0, 0) circle [x radius= 3.35, y radius= 3.35]   ;
\draw [color={rgb, 255:red, 0; green, 0; blue, 0 }  ,draw opacity=1 ][line width=0.75]    (91.5,151.5) ;
\draw [shift={(91.5,151.5)}, rotate = 0] [color={rgb, 255:red, 0; green, 0; blue, 0 }  ,draw opacity=1 ][fill={rgb, 255:red, 0; green, 0; blue, 0 }  ,fill opacity=1 ][line width=0.75]      (0, 0) circle [x radius= 3.35, y radius= 3.35]   ;
\draw [shift={(91.5,151.5)}, rotate = 0] [color={rgb, 255:red, 0; green, 0; blue, 0 }  ,draw opacity=1 ][line width=0.75]      (0, 0) circle [x radius= 3.35, y radius= 3.35]   ;
\draw [color={rgb, 255:red, 208; green, 2; blue, 27 }  ,draw opacity=1 ][line width=0.75]    (440.13,150.16) ;
\draw [shift={(440.13,150.16)}, rotate = 0] [color={rgb, 255:red, 208; green, 2; blue, 27 }  ,draw opacity=1 ][fill={rgb, 255:red, 208; green, 2; blue, 27 }  ,fill opacity=1 ][line width=0.75]      (0, 0) circle [x radius= 3.35, y radius= 3.35]   ;
\draw [shift={(440.13,150.16)}, rotate = 0] [color={rgb, 255:red, 208; green, 2; blue, 27 }  ,draw opacity=1 ][line width=0.75]      (0, 0) circle [x radius= 3.35, y radius= 3.35]   ;
\draw [color={rgb, 255:red, 208; green, 2; blue, 27 }  ,draw opacity=1 ][line width=0.75]    (381.33,150.36) ;
\draw [shift={(381.33,150.36)}, rotate = 0] [color={rgb, 255:red, 208; green, 2; blue, 27 }  ,draw opacity=1 ][fill={rgb, 255:red, 208; green, 2; blue, 27 }  ,fill opacity=1 ][line width=0.75]      (0, 0) circle [x radius= 2.68, y radius= 2.68]   ;
\draw [shift={(381.33,150.36)}, rotate = 0] [color={rgb, 255:red, 208; green, 2; blue, 27 }  ,draw opacity=1 ][line width=0.75]      (0, 0) circle [x radius= 2.68, y radius= 2.68]   ;
\draw [color={rgb, 255:red, 74; green, 144; blue, 226 }  ,draw opacity=1 ][line width=0.75]    (355.99,45.37) ;
\draw [shift={(355.99,45.37)}, rotate = 0] [color={rgb, 255:red, 74; green, 144; blue, 226 }  ,draw opacity=1 ][fill={rgb, 255:red, 74; green, 144; blue, 226 }  ,fill opacity=1 ][line width=0.75]      (0, 0) circle [x radius= 3.69, y radius= 3.69]   ;
\draw [shift={(355.99,45.37)}, rotate = 0] [color={rgb, 255:red, 74; green, 144; blue, 226 }  ,draw opacity=1 ][line width=0.75]      (0, 0) circle [x radius= 3.69, y radius= 3.69]   ;
\draw [color={rgb, 255:red, 74; green, 144; blue, 226 }  ,draw opacity=1 ][line width=0.75]    (105,114.23) ;
\draw [shift={(105,114.23)}, rotate = 0] [color={rgb, 255:red, 74; green, 144; blue, 226 }  ,draw opacity=1 ][fill={rgb, 255:red, 74; green, 144; blue, 226 }  ,fill opacity=1 ][line width=0.75]      (0, 0) circle [x radius= 3.69, y radius= 3.69]   ;
\draw [shift={(105,114.23)}, rotate = 0] [color={rgb, 255:red, 74; green, 144; blue, 226 }  ,draw opacity=1 ][line width=0.75]      (0, 0) circle [x radius= 3.69, y radius= 3.69]   ;
\draw [color={rgb, 255:red, 74; green, 144; blue, 226 }  ,draw opacity=1 ][line width=0.75]    (393,150.36) ;
\draw [shift={(393,150.36)}, rotate = 0] [color={rgb, 255:red, 74; green, 144; blue, 226 }  ,draw opacity=1 ][fill={rgb, 255:red, 74; green, 144; blue, 226 }  ,fill opacity=1 ][line width=0.75]      (0, 0) circle [x radius= 3.69, y radius= 3.69]   ;
\draw [shift={(393,150.36)}, rotate = 0] [color={rgb, 255:red, 74; green, 144; blue, 226 }  ,draw opacity=1 ][line width=0.75]      (0, 0) circle [x radius= 3.69, y radius= 3.69]   ;
\draw [color={rgb, 255:red, 74; green, 144; blue, 226 }  ,draw opacity=1 ][line width=0.75]    (451,167) ;
\draw [shift={(451,167)}, rotate = 0] [color={rgb, 255:red, 74; green, 144; blue, 226 }  ,draw opacity=1 ][fill={rgb, 255:red, 74; green, 144; blue, 226 }  ,fill opacity=1 ][line width=0.75]      (0, 0) circle [x radius= 3.69, y radius= 3.69]   ;
\draw [shift={(451,167)}, rotate = 0] [color={rgb, 255:red, 74; green, 144; blue, 226 }  ,draw opacity=1 ][line width=0.75]      (0, 0) circle [x radius= 3.69, y radius= 3.69]   ;

\draw [color={rgb, 255:red, 74; green, 144; blue, 226 }  ,draw opacity=1 ][line width=0.75]    (550,190.23) ;
\draw [shift={(550,190.23)}, rotate = 0] [color={rgb, 255:red, 74; green, 144; blue, 226 }  ,draw opacity=1 ][fill={rgb, 255:red, 74; green, 144; blue, 226 }  ,fill opacity=1 ][line width=0.75]      (0, 0) circle [x radius= 3.69, y radius= 3.69]   ;
\draw [shift={(550,190.23)}, rotate = 0] [color={rgb, 255:red, 74; green, 144; blue, 226 }  ,draw opacity=1 ][line width=0.75]      (0, 0) circle [x radius= 3.69, y radius= 3.69]   ;

\draw    (296,-1.9) -- (331,151.37) ;
\draw    (366,-1.9) -- (331,151.37) ;

\draw (325,164.63) node [anchor=north west][inner sep=0.75pt]    {$0$};
\draw (369,11.4) node [anchor=north west][inner sep=0.75pt]    {$\mathbb{E}(0 ,w ,\varepsilon )$};
\draw (74,163.4) node [anchor=north west][inner sep=0.75pt]    {$-Ru$};
\draw (550,160.4) node [anchor=north west][inner sep=0.75pt]    {$Ru$};
\draw (351,54.4) node [anchor=north west][inner sep=0.75pt]    {$\textcolor[rgb]{0.29,0.56,0.89}{q}\textcolor[rgb]{0.29,0.56,0.89}{^{++}}$};
\draw (101,77.4) node [anchor=north west][inner sep=0.75pt]    {$\textcolor[rgb]{0.29,0.56,0.89}{q}\textcolor[rgb]{0.29,0.56,0.89}{^{-}}$};
\draw (442.13,173.03) node [anchor=north west][inner sep=0.75pt]    {$\textcolor[rgb]{0.29,0.56,0.89}{q}\textcolor[rgb]{0.29,0.56,0.89}{_{1}^{+}}$};
\draw (383.33,153.76) node [anchor=north west][inner sep=0.75pt]    {$\textcolor[rgb]{0.29,0.56,0.89}{q}\textcolor[rgb]{0.29,0.56,0.89}{_{2}^{+}}$};
\draw (388,255.4) node [anchor=north west][inner sep=0.75pt]    {$\textcolor[rgb]{0.82,0.01,0.11}{B(0,R_1)}$};
\draw (538,205.4) node [anchor=north west][inner sep=0.75pt]    {$\textcolor[rgb]{0.29,0.56,0.89}{q}\textcolor[rgb]{0.29,0.56,0.89}{^{+}}$};
\draw (528,85.4) node [anchor=north west][inner sep=0.75pt]    {$B(Ru,\lambda R\varepsilon_R)$};
\draw (28,205.4) node [anchor=north west][inner sep=0.75pt]    {$B(-Ru,\lambda R\varepsilon_{-R})$};

\end{tikzpicture}
\selectlanguage{english}
\caption{Case~\eqref{dwa}}
\label{f2}
\end{figure}

\vspace{-1.5cm}
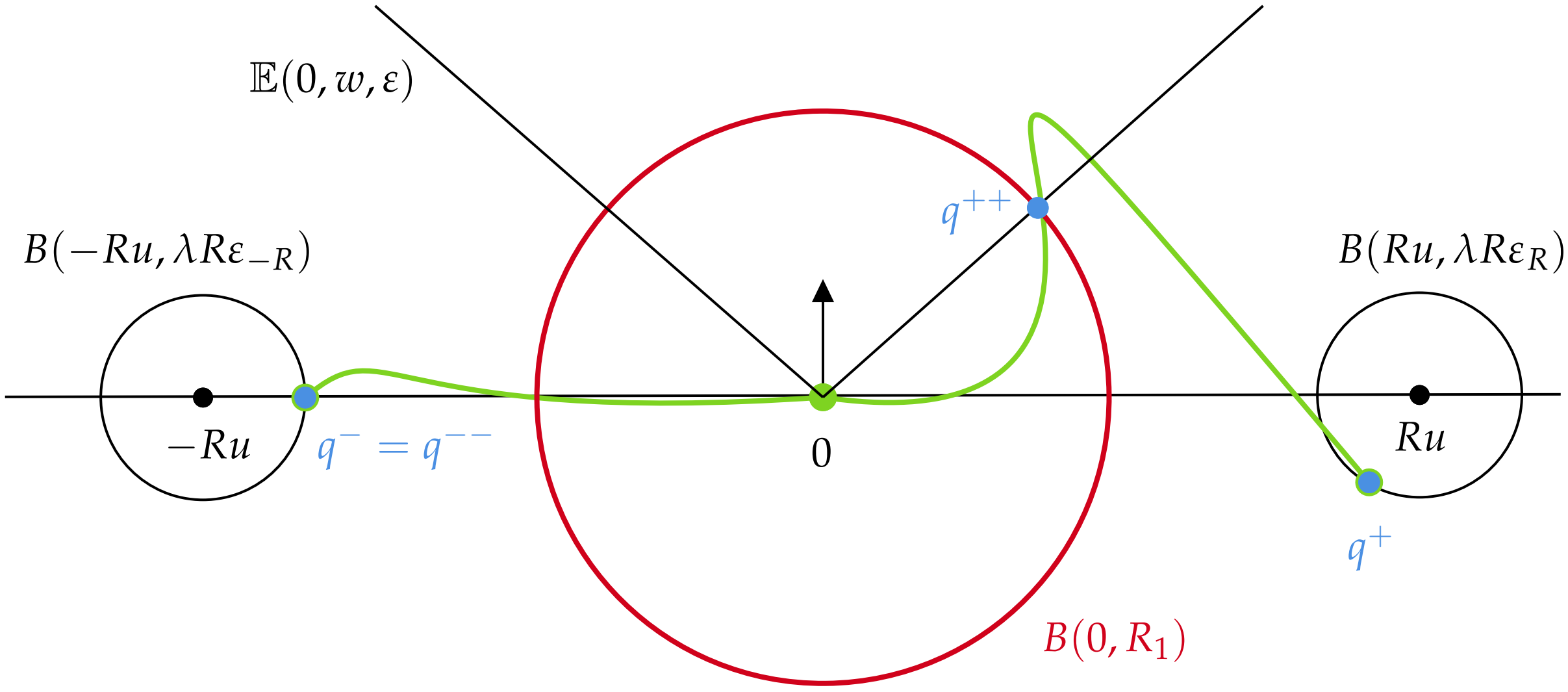
\begin{figure}[H]
\centering
\tikzset{every picture/.style={line width=0.75pt}}

\begin{tikzpicture}[x=0.75pt,y=0.75pt,yscale=-1,xscale=1]

\draw    (15,151.23) -- (616,150.23) ;
\draw    (331,108.73) -- (331,151.37) ;
\draw [shift={(331,105.73)}, rotate = 90] [fill={rgb, 255:red, 0; green, 0; blue, 0 }  ][line width=0.08]  [draw opacity=0] (8.93,-4.29) -- (0,0) -- (8.93,4.29) -- cycle    ;
\draw   (522,150.5) .. controls (522,128.68) and (539.68,111) .. (561.5,111) .. controls (583.32,111) and (601,128.68) .. (601,150.5) .. controls (601,172.32) and (583.32,190) .. (561.5,190) .. controls (539.68,190) and (522,172.32) .. (522,150.5) -- cycle ;
\draw   (52,151.5) .. controls (52,129.68) and (69.68,112) .. (91.5,112) .. controls (113.32,112) and (131,129.68) .. (131,151.5) .. controls (131,173.32) and (113.32,191) .. (91.5,191) .. controls (69.68,191) and (52,173.32) .. (52,151.5) -- cycle ;
\draw [color={rgb, 255:red, 126; green, 211; blue, 33 }  ,draw opacity=1 ][line width=1.5]    (131,151.5) .. controls (164,123.23) and (158,163.23) .. (331,151.37) ;
\draw [shift={(331,151.37)}, rotate = 356.08] [color={rgb, 255:red, 126; green, 211; blue, 33 }  ,draw opacity=1 ][fill={rgb, 255:red, 126; green, 211; blue, 33 }  ,fill opacity=1 ][line width=1.5]      (0, 0) circle [x radius= 4.36, y radius= 4.36]   ;
\draw [shift={(131,151.5)}, rotate = 319.42] [color={rgb, 255:red, 126; green, 211; blue, 33 }  ,draw opacity=1 ][fill={rgb, 255:red, 126; green, 211; blue, 33 }  ,fill opacity=1 ][line width=1.5]      (0, 0) circle [x radius= 4.36, y radius= 4.36]   ;
\draw [color={rgb, 255:red, 126; green, 211; blue, 33 }  ,draw opacity=1 ][line width=1.5]    (331,151.37) .. controls (534,181.23) and (276,-130.77) .. (542,184.23) ;
\draw [shift={(542,184.23)}, rotate = 49.82] [color={rgb, 255:red, 126; green, 211; blue, 33 }  ,draw opacity=1 ][fill={rgb, 255:red, 126; green, 211; blue, 33 }  ,fill opacity=1 ][line width=1.5]      (0, 0) circle [x radius= 4.36, y radius= 4.36]   ;
\draw [shift={(331,151.37)}, rotate = 8.37] [color={rgb, 255:red, 126; green, 211; blue, 33 }  ,draw opacity=1 ][fill={rgb, 255:red, 126; green, 211; blue, 33 }  ,fill opacity=1 ][line width=1.5]      (0, 0) circle [x radius= 4.36, y radius= 4.36]   ;
\draw  [color={rgb, 255:red, 208; green, 2; blue, 27 }  ,draw opacity=1 ][line width=1.5]  (220.48,152.45) .. controls (219.89,91.41) and (268.88,41.45) .. (329.92,40.85) .. controls (390.96,40.25) and (440.92,89.25) .. (441.52,150.29) .. controls (442.11,211.32) and (393.12,261.29) .. (332.08,261.88) .. controls (271.04,262.48) and (221.08,213.48) .. (220.48,152.45) -- cycle ;
\draw [color={rgb, 255:red, 0; green, 0; blue, 0 }  ,draw opacity=1 ][line width=0.75]    (561.5,150.5) ;
\draw [shift={(561.5,150.5)}, rotate = 0] [color={rgb, 255:red, 0; green, 0; blue, 0 }  ,draw opacity=1 ][fill={rgb, 255:red, 0; green, 0; blue, 0 }  ,fill opacity=1 ][line width=0.75]      (0, 0) circle [x radius= 3.35, y radius= 3.35]   ;
\draw [shift={(561.5,150.5)}, rotate = 0] [color={rgb, 255:red, 0; green, 0; blue, 0 }  ,draw opacity=1 ][line width=0.75]      (0, 0) circle [x radius= 3.35, y radius= 3.35]   ;
\draw [color={rgb, 255:red, 0; green, 0; blue, 0 }  ,draw opacity=1 ][line width=0.75]    (91.5,151.5) ;
\draw [shift={(91.5,151.5)}, rotate = 0] [color={rgb, 255:red, 0; green, 0; blue, 0 }  ,draw opacity=1 ][fill={rgb, 255:red, 0; green, 0; blue, 0 }  ,fill opacity=1 ][line width=0.75]      (0, 0) circle [x radius= 3.35, y radius= 3.35]   ;
\draw [shift={(91.5,151.5)}, rotate = 0] [color={rgb, 255:red, 0; green, 0; blue, 0 }  ,draw opacity=1 ][line width=0.75]      (0, 0) circle [x radius= 3.35, y radius= 3.35]   ;
\draw [color={rgb, 255:red, 74; green, 144; blue, 226 }  ,draw opacity=1 ][line width=0.75]    (131,151.5) ;
\draw [shift={(131,151.5)}, rotate = 0] [color={rgb, 255:red, 74; green, 144; blue, 226 }  ,draw opacity=1 ][fill={rgb, 255:red, 74; green, 144; blue, 226 }  ,fill opacity=1 ][line width=0.75]      (0, 0) circle [x radius= 3.69, y radius= 3.69]   ;
\draw [shift={(131,151.5)}, rotate = 0] [color={rgb, 255:red, 74; green, 144; blue, 226 }  ,draw opacity=1 ][line width=0.75]      (0, 0) circle [x radius= 3.69, y radius= 3.69]   ;
\draw    (158,0.23) -- (331,151.37) ;
\draw    (501,0.23) -- (331,151.37) ;
\draw [color={rgb, 255:red, 74; green, 144; blue, 226 }  ,draw opacity=1 ][line width=0.75]    (414,78.23) ;
\draw [shift={(414,78.23)}, rotate = 0] [color={rgb, 255:red, 74; green, 144; blue, 226 }  ,draw opacity=1 ][fill={rgb, 255:red, 74; green, 144; blue, 226 }  ,fill opacity=1 ][line width=0.75]      (0, 0) circle [x radius= 3.69, y radius= 3.69]   ;
\draw [shift={(414,78.23)}, rotate = 0] [color={rgb, 255:red, 74; green, 144; blue, 226 }  ,draw opacity=1 ][line width=0.75]      (0, 0) circle [x radius= 3.69, y radius= 3.69]   ;
\draw [color={rgb, 255:red, 74; green, 144; blue, 226 }  ,draw opacity=1 ][line width=0.75]    (220.48,152.45) ;
\draw [color={rgb, 255:red, 74; green, 144; blue, 226 }  ,draw opacity=1 ][line width=0.75]    (542,184.23) ;
\draw [shift={(542,184.23)}, rotate = 0] [color={rgb, 255:red, 74; green, 144; blue, 226 }  ,draw opacity=1 ][fill={rgb, 255:red, 74; green, 144; blue, 226 }  ,fill opacity=1 ][line width=0.75]      (0, 0) circle [x radius= 3.69, y radius= 3.69]   ;
\draw [shift={(542,184.23)}, rotate = 0] [color={rgb, 255:red, 74; green, 144; blue, 226 }  ,draw opacity=1 ][line width=0.75]      (0, 0) circle [x radius= 3.69, y radius= 3.69]   ;

\draw (325,165.63) node [anchor=north west][inner sep=0.75pt]    {$0$};
\draw (108,19.4) node [anchor=north west][inner sep=0.75pt]    {$\mathbb{E}(0 ,w ,\varepsilon )$};
\draw (75,162.4) node [anchor=north west][inner sep=0.75pt]    {$-Ru$};
\draw (550,159.4) node [anchor=north west][inner sep=0.75pt]    {$Ru$};
\draw (375,69.63) node [anchor=north west][inner sep=0.75pt]    {$\textcolor[rgb]{0.29,0.56,0.89}{q^{++}}$};
\draw (134,160.4) node [anchor=north west][inner sep=0.75pt]    {$\textcolor[rgb]{0.29,0.56,0.89}{q^{-}=q^{--}}$};
\draw (532,199.4) node [anchor=north west][inner sep=0.75pt]    {$\textcolor[rgb]{0.29,0.56,0.89}{q^{+}}$};
\draw (414,235.4) node [anchor=north west][inner sep=0.75pt]    {$\textcolor[rgb]{0.82,0.01,0.11}{B(0,R_1)}$};
\draw (528,85.4) node [anchor=north west][inner sep=0.75pt]    {$B(Ru,\lambda R\varepsilon_R)$};
\draw (20,85.4) node [anchor=north west][inner sep=0.75pt]    {$B(-Ru,\lambda R\varepsilon_{-R})$};

\end{tikzpicture}

\selectlanguage{english}
\caption{Case~\eqref{trzy}: a configuration in which the arcs avoid the cone near $0$. 
}
\label{f3}
\end{figure}

%
    We will show that, in each case, the upper Hausdorff density of $g(B)$ at the origin exceeds $1$, contradicting~\eqref{zaloz}.

    In the next step, we address the first two cases.
    \needspace{5\baselineskip}
    \begin{center}
    \textbf{Step 2. Cases \eqref{jeden} and \eqref{dwa}}
    \end{center}
    We start with Case~\eqref{jeden}.

    In this case, the subarcs $[0,q^{++}]$ and $[0,q^{--}]$ are contained in $\overline{\mathbb E(0,w,\varepsilon)}$. Put $a\coloneqq\min\{\|q^{++}\|,\|q^{--}\|\}>0$ and take any $0<R_1\le\min\{a,R\}$. The initial part of $[0,q^{++}]$ ending at its first intersection with $\partial B(0,R_1)$ is contained in $\overline{B(0,R_1)}$ and has Hausdorff measure at least $R_1$, by Proposition~\ref{haus}.

    Set
    \[
    R_2\coloneqq\frac12R_1,
    \quad
    R_4\coloneqq\frac14R_1.
    \]
    For \(j\in\{1,2,4\}\), put
    \[
    q_j^+\coloneqq g(R_jv),
    \quad
    q_j^-\coloneqq g(-R_jv).
    \]
    By \eqref{waznee} we have (see Figure~\ref{f1})
    \[
    q_j^+\in
    g(B)\cap
    \partial B(R_ju,R_j\lambda\varepsilon_{R_j}),\quad
    q_j^-\in
    g(B)\cap
    \partial B(-R_ju,R_j\lambda\varepsilon_{-R_j}).
    \]
    We connect these points using the images of the intervening parameter intervals under $t\mapsto g(tv)$ and $t\mapsto g(-tv)$. The green curves in Figures~\ref{f1} and~\ref{f2} indicate possible
    connecting arcs.

    Follow the positive path from $q_4^+$ towards $q_1^+$, via $q_2^+$, stopping at its first intersection with $\partial B(0,R_1)$ if one occurs. Otherwise, keep the whole path up to $q_1^+$. 
    The resulting image is connected and lies in $\overline{B(0,R_1)}$. For $0<t\le R_2$, by \eqref{taylor} and the triangle inequality we then have
    \[
    \|g(tv)\|=\|tu+R_0(tv)\|\le t(1+\lambda\varepsilon_t)
       <\frac54t\le\frac 54R_2=\frac58R_1<R_1,
    \]
    so the retained path passes through $q_2^+$ before it can meet the sphere. If it ends at $q_1^+$, the distance between its endpoints is at least
    \[
    \begin{aligned}
    \|R_1u-R_4u\|-R_1\lambda\varepsilon_{R_1}
                         -R_4\lambda\varepsilon_{R_4}
    &=\frac34R_1-R_1\lambda\varepsilon_{R_1}
                         -\frac14R_1\lambda\varepsilon_{R_4}.
    \end{aligned}
    \]
    If it is stopped on $\partial B(0,R_1)$, the distance from $q_4^+$ to its endpoint is at least
    \[
    R_1-\|q_4^+\|
       \ge\frac34R_1-\frac14R_1\lambda\varepsilon_{R_4},
    \]
    which is an even stronger bound. Proposition~\ref{haus} therefore gives the first lower bound for the measure of the connected piece in either situation. The same argument applies to the negative path from $q_4^-$ through $q_2^-$ towards $q_1^-$. If this path reaches $q_1^-$ before meeting $\partial B(0,R_1)$, the distance between its endpoints is at least
    \[
    \|{-R_1u}-(-R_4u)\|
    -
    R_1\lambda\varepsilon_{-R_1}
    -
    R_4\lambda\varepsilon_{-R_4}
    =
    \frac34R_1
    -
    R_1\lambda\varepsilon_{-R_1}
    -
    \frac14R_1\lambda\varepsilon_{-R_4}.
    \]
    If it is stopped on $\partial B(0,R_1)$, its diameter is at least
    \[
    R_1-\|q_4^-\|
    \geq
    \frac34R_1
    -
    \frac14R_1\lambda\varepsilon_{-R_4},
    \]
    which is again an even stronger bound.
    
    The positive and negative pieces lie in the disjoint cones about $u$ and $-u$. They are also disjoint from the part of $[0,q^{++}]$ contained in the closed cone about $w$. Hence their Hausdorff measures can be added. Therefore
    \[
    \begin{split}
    \hau^1\left(g(B)\cap\overline{B(0,R_1)}\right)
    \geq {}&
    R_1+
    \left(
    \frac34R_1-
    R_1\lambda\varepsilon_{R_1}
    -
    \frac14R_1\lambda\varepsilon_{R_4}
    \right)\\
    &+
    \left(
    \frac34R_1
    -
    R_1\lambda\varepsilon_{-R_1}
    -
    \frac14R_1\lambda\varepsilon_{-R_4}
    \right).
    \end{split}
    \]
    Consequently,
    \[
    \begin{split}
    \theta_1^*(g(B),0)
    &\geq
    \limsup_{R_1\to0^+}
    \left[
    \frac54
    -
    \frac{\lambda}{2}
    \left(
    \varepsilon_{R_1}
    +
    \varepsilon_{-R_1}
    \right)
    -
    \frac{\lambda}{8}
    \left(
    \varepsilon_{R_4}
    +
    \varepsilon_{-R_4}
    \right)
    \right]=
    \frac54>1.
    \end{split}
    \]
    Indeed, $R_4=R_1/4$, and
    \[
    \varepsilon_{R_1},
    \varepsilon_{-R_1},
    \varepsilon_{R_4},
    \varepsilon_{-R_4}
    \xrightarrow{R_1\to 0^+}0.
    \]
    This contradicts \eqref{zaloz}.

    \vspace{0.2cm}
    Now we address Case~\eqref{dwa}. Put $a\coloneqq\|q^{++}\|$, $b\coloneqq\|q^-\|$, and take any $0<R_1\le\min\{R,a,b\}$. As before, the initial part of $[0,q^{++}]$ up to its first intersection with $\partial B(0,R_1)$ lies in the closed cone about $w$ and has measure at least $R_1$.

    For the negative contribution, use the path $t\mapsto g(-tv)$, $0\le t\le R$, from $0$ to $q^-$. Its nonzero points lie in $\mathbb E(0,-u,\delta)$. Since $\|q^-\|\ge R_1$, its initial part up to the first intersection with $\partial B(0,R_1)$ is connected, lies in the closed ball, and has measure at least $R_1$. This path need not coincide with the arbitrary arc $[0,q^-]$ used to distinguish the cases. The two contributions are disjoint apart from the point $0$, and hence
    \begin{equation}\label{proste}
    \hau^1\left(g(B)\cap\overline{B(0,R_1)}\right)\ge R_1+R_1.
    \end{equation}

    To obtain the additional contribution shown on the right of Figure~\ref{f2}, put $R_2\coloneqq R_1/2$ and $q_1^+\coloneqq g(R_1v)$, $q_2^+\coloneqq g(R_2v)$.
    Follow $t\mapsto g(tv)$ from $q_2^+$ towards $q_1^+$, using the same stopping rule as in Case~\eqref{jeden}. The retained image is connected and lies in $\overline{B(0,R_1)}\cap\mathbb E(0,u,\delta)$. If it ends at $q_1^+$, its diameter is at least
    \[
    \frac12R_1-R_1\lambda\varepsilon_{R_1}
                     -\frac12R_1\lambda\varepsilon_{R_2}.
    \]
    If it ends on the sphere, its diameter is at least $R_1-\|q_2^+\|$, which is no smaller than this bound. This third piece is disjoint from the two pieces already counted, because it lies in the positive cone. Consequently,
    \begin{equation*}
    \begin{split}
    \hau^1\left(g(B)\cap\overline{B(0,R_1)}\right)
    &\ge R_1+R_1+\frac12R_1-R_1\lambda\varepsilon_{R_1}
                            -\frac12R_1\lambda\varepsilon_{R_2},\\
    \theta_1^*(g(B),0)
    &\ge\limsup_{R_1\to0^+}
       \left(\frac54-\frac\lambda2\varepsilon_{R_1}
                        -\frac\lambda4\varepsilon_{R_2}\right)
     =\frac54>1.
    \end{split}
    \end{equation*}
    This again contradicts~\eqref{zaloz}.

    It remains to consider Case~\eqref{trzy}. Here an initial subarc inside the cone about $w$ is not available. We instead use the assumption $w\in\operatorname{Tan}(g(B),0)$ to find points of $g(B)$ tending to $0$ in direction $w$, and connect successive points by arcs. Their parts inside the cone give an additional contribution to the measure. We combine this with the contribution outside the cone at the same radii.

    \needspace{5\baselineskip}
    \begin{center}
    \textbf{Step 3. Case~\eqref{trzy}}
    \end{center}
    Set \(F\coloneqq\mathbb R^2\setminus\mathbb E(0,w,\varepsilon)\). Recall that \(0<\varepsilon<1/2\). We begin with a lower bound outside the cone, valid at every sufficiently small radius.

    Take any $0<R_1\le\min\{\|q^+\|,\|q^-\|\}$. Follow the two paths $t\mapsto g(tv)$ and $t\mapsto g(-tv)$, $0\le t\le R$, from $0$ until each first meets $\partial B(0,R_1)$. Each retained image is connected, lies in $\overline{B(0,R_1)}$, and has diameter at least $R_1$. Their nonzero points lie in the disjoint cones about $u$ and $-u$, both contained in $F$. The two images can intersect only at $0$, a set of zero $\hau^1$ measure. Proposition~\ref{haus} therefore gives
    \begin{equation}\label{naef}
    \hau^1\left((g(B)\cap F)\cap\overline{B(0,R_1)}\right)
       \ge R_1+R_1=2R_1.
    \end{equation}
    In Figure~\ref{f3}, $q^{++}$ and $q^{--}$ label the endpoints of the cone-avoiding subarcs in the illustrated configuration. The estimate does not require the arbitrary arcs used in the case distinction to avoid the cone near $0$.

    In particular,
    \begin{equation}\label{naeff}
    \begin{split}
    \theta_1^*(g(B)\cap F,0)
    &=\limsup_{R_1\to0^+}
       \frac{\hau^1((g(B)\cap F)\cap\overline{B(0,R_1)})}{2R_1}\ge1.
    \end{split}
    \end{equation}
    Next we show how to obtain a positive contribution inside the cone. We will then add the measure estimates at the same radii using~\eqref{naef} and reach a contradiction.

    By the definition of a tangent vector, $w\in\text{Tan}(g(B),0)$ means that for every $t>0$ there exist $\tau>0$ and $z\in g(B)$ such that $\|z\|< t$ and $\|\tau z-w\|<t$. In particular (by the definition of a directional cone) for every $t\in(0,\varepsilon)$ there exists a point $z$ such that $z\in g(B)\cap B(0,t)\cap\mathbb{E}(0,w,\varepsilon)$.
    We first consider a special case and then the general case, both follow the same geometric idea.

     \Needspace{6\baselineskip}
    \medskip\noindent\textbf{Situation 1. Points on the ray in direction $w$.}

     Assume first that
     \[\forall_{t\in(0,\varepsilon)} \exists_{\alpha\in(0,t)} \text{ such that }
     z=\alpha w\in g(B)\cap B(0,t)\cap\mathbb{E}(0,w,\varepsilon).\]
     This means that there exists a sequence $(\alpha_n)_{n=1}^{\infty}$ of positive numbers bounded by $\varepsilon$ such that $\alpha_n\downarrow 0$ and for every $n\in\nat$
     \[z_n\coloneqq\alpha_n w\in g(B)\cap \mathbb{E}(0,w,\varepsilon).\]
     Since $\alpha_n$ converges to $0$ we may assume that $\alpha_{n+1}\leq\frac{1}{2}\alpha_n$ for all $n\in\nat$ after passing to a subsequence (denoted the same), if necessary.
%
     Since $g(B)$ is arcwise connected, we can connect $z_n$ and $z_{n+1}$ by an arc in $g(B)$. We estimate its measure inside $B(0,2\alpha_n)\cap\mathbb E(0,w,\varepsilon)$. The radius is $2\alpha_n$, as measuring only in $B(0,\alpha_n)$ would put the starting point on the sphere, allowing the possibility of the arc immediately leaving the measuring ball.

     Follow the arc from $z_n$ until it first reaches $z_{n+1}$, the boundary of the cone, or the sphere $\partial B(0,2\alpha_n)$. The retained part is connected. If it reaches $z_{n+1}$ first, the distance between its endpoints is $\alpha_n-\alpha_{n+1}\ge\alpha_n/2$. If it reaches the boundary of the cone first, that distance is at least $\varepsilon\alpha_n$, the perpendicular distance shown in Figure~\ref{f4}. Here the half-angle is $\arcsin\varepsilon$ and $\|w\|=1$.
     If it reaches the outer sphere first, the distance is at least $2\alpha_n-\|z_n\|=\alpha_n$. Since $\varepsilon<1/2$, Proposition~\ref{haus} gives a lower bound of $\varepsilon\alpha_n$ in all three alternatives. Thus
     \[
     \hau^1\left(g(B)\cap\overline{B(0,2\alpha_n)}
                    \cap\mathbb E(0,w,\varepsilon)\right)
       \ge\varepsilon\alpha_n.
     \]
     For all sufficiently large $n$, we may apply~\eqref{naef} at $R_1=2\alpha_n$. The parts inside and outside the cone are disjoint, so
     \[
     \hau^1\left(g(B)\cap\overline{B(0,2\alpha_n)}\right)
       \ge4\alpha_n+\varepsilon\alpha_n.
     \]
     Consequently,
     \[
     \theta_1^*(g(B),0)
       \ge\limsup_{n\to\infty}
          \frac{4\alpha_n+\varepsilon\alpha_n}{4\alpha_n}
       =1+\frac{\varepsilon}{4}>1,
     \]
     which contradicts~\eqref{zaloz}.
\vspace{0.1cm}
\begin{figure}[H]
\centering

\tikzset{every picture/.style={line width=0.75pt}}

\begin{tikzpicture}[x=0.75pt,y=0.75pt,yscale=-1,xscale=1]

\draw    (331,252.23) -- (329.99,7.85) ;
\draw [shift={(329.99,5.85)}, rotate = 89.76] [color={rgb, 255:red, 0; green, 0; blue, 0 }  ][line width=0.75]    (10.93,-3.29) .. controls (6.95,-1.4) and (3.31,-0.3) .. (0,0) .. controls (3.31,0.3) and (6.95,1.4) .. (10.93,3.29)   ;
\draw    (193,0.23) -- (331,252.23) ;
\draw [color={rgb, 255:red, 144; green, 19; blue, 254 }  ,draw opacity=1 ][line width=1.5]    (330,30.23) -- (294.58,49.25) -- (283.13,55.41) -- (236.91,80.22) ;
\draw [color={rgb, 255:red, 0; green, 0; blue, 0 }  ,draw opacity=1 ][line width=0.75]    (331,252.23) ;
\draw [shift={(331,252.23)}, rotate = 0] [color={rgb, 255:red, 0; green, 0; blue, 0 }  ,draw opacity=1 ][fill={rgb, 255:red, 0; green, 0; blue, 0 }  ,fill opacity=1 ][line width=0.75]      (0, 0) circle [x radius= 3.35, y radius= 3.35]   ;
\draw [shift={(331,252.23)}, rotate = 0] [color={rgb, 255:red, 0; green, 0; blue, 0 }  ,draw opacity=1 ][line width=0.75]      (0, 0) circle [x radius= 3.35, y radius= 3.35]   ;
\draw [color={rgb, 255:red, 0; green, 0; blue, 0 }  ,draw opacity=1 ][line width=0.75]    (329.99,5.85) ;
\draw [shift={(329.99,5.85)}, rotate = 0] [color={rgb, 255:red, 0; green, 0; blue, 0 }  ,draw opacity=1 ][fill={rgb, 255:red, 0; green, 0; blue, 0 }  ,fill opacity=1 ][line width=0.75]      (0, 0) circle [x radius= 3.35, y radius= 3.35]   ;
\draw [shift={(329.99,5.85)}, rotate = 0] [color={rgb, 255:red, 0; green, 0; blue, 0 }  ,draw opacity=1 ][line width=0.75]      (0, 0) circle [x radius= 3.35, y radius= 3.35]   ;
\draw [color={rgb, 255:red, 0; green, 0; blue, 0 }  ,draw opacity=1 ][line width=0.75]    (330,30.23) ;
\draw [shift={(330,30.23)}, rotate = 0] [color={rgb, 255:red, 0; green, 0; blue, 0 }  ,draw opacity=1 ][fill={rgb, 255:red, 0; green, 0; blue, 0 }  ,fill opacity=1 ][line width=0.75]      (0, 0) circle [x radius= 3.35, y radius= 3.35]   ;
\draw [shift={(330,30.23)}, rotate = 0] [color={rgb, 255:red, 0; green, 0; blue, 0 }  ,draw opacity=1 ][line width=0.75]      (0, 0) circle [x radius= 3.35, y radius= 3.35]   ;
\draw [color={rgb, 255:red, 0; green, 0; blue, 0 }  ,draw opacity=1 ][line width=0.75]    (331.28,161.76) ;
\draw [shift={(331.28,161.76)}, rotate = 0] [color={rgb, 255:red, 0; green, 0; blue, 0 }  ,draw opacity=1 ][fill={rgb, 255:red, 0; green, 0; blue, 0 }  ,fill opacity=1 ][line width=0.75]      (0, 0) circle [x radius= 3.35, y radius= 3.35]   ;
\draw [shift={(331.28,161.76)}, rotate = 0] [color={rgb, 255:red, 0; green, 0; blue, 0 }  ,draw opacity=1 ][line width=0.75]      (0, 0) circle [x radius= 3.35, y radius= 3.35]   ;
\draw [color={rgb, 255:red, 144; green, 19; blue, 254 }  ,draw opacity=1 ][line width=1.5]    (292.79,181.95) -- (331.28,161.76) ;
\draw   (303.77,176.26) -- (311,190.23) -- (300.02,195.92) ;
\draw    (227.61,62.25) -- (329.99,5.85) ;
\draw   (239.83,55.93) -- (244.83,65.59) -- (232,72.23) ;
\draw   (249.13,73.9) -- (254.91,85.07) -- (243,91.23) ;
\draw    (29.88,252.68) -- (632.13,251.78) ;
\draw    (472,-0.77) -- (331,252.23) ;
\draw [color={rgb, 255:red, 245; green, 166; blue, 35 }  ,draw opacity=1 ][line width=1.5]    (330,30.23) -- (331.28,161.76) ;

\draw (339.99,9.25) node [anchor=north west][inner sep=0.75pt]    {$\mathnormal{w}$};
\draw (332,33.63) node [anchor=north west][inner sep=0.75pt]    {$\alpha _{n} w$};
\draw (332.49,132.44) node [anchor=north west][inner sep=0.75pt]    {$\alpha_{n+1} w$};
\draw (264,15.4) node [anchor=north west][inner sep=0.75pt]    {$\varepsilon$};
\draw (469,16.63) node [anchor=north west][inner sep=0.75pt]    {$\mathbb{E}( 0,w,\varepsilon)$};
\draw (333,80.4) node [anchor=north west][inner sep=0.75pt]    {$\mathnormal{\textcolor[rgb]{0.96,0.65,0.14}{\alpha _{n} -\alpha_{n+1}}}$};
\draw (276,63.4) node [anchor=north west][inner sep=0.75pt]    {$\mathnormal{\textcolor[rgb]{0.56,0.07,1}{\varepsilon}\textcolor[rgb]{0.56,0.07,1}{\alpha }\textcolor[rgb]{0.56,0.07,1}{_{n}}}$};
\draw (283,143.4) node [anchor=north west][inner sep=0.75pt]    {$\mathnormal{\textcolor[rgb]{0.56,0.07,1}{\varepsilon}\textcolor[rgb]{0.56,0.07,1}{\alpha }\textcolor[rgb]{0.56,0.07,1}{_{n+1}}}$};

\end{tikzpicture}

\selectlanguage{english}
\caption{Situation 1}
\label{f4}
\end{figure}
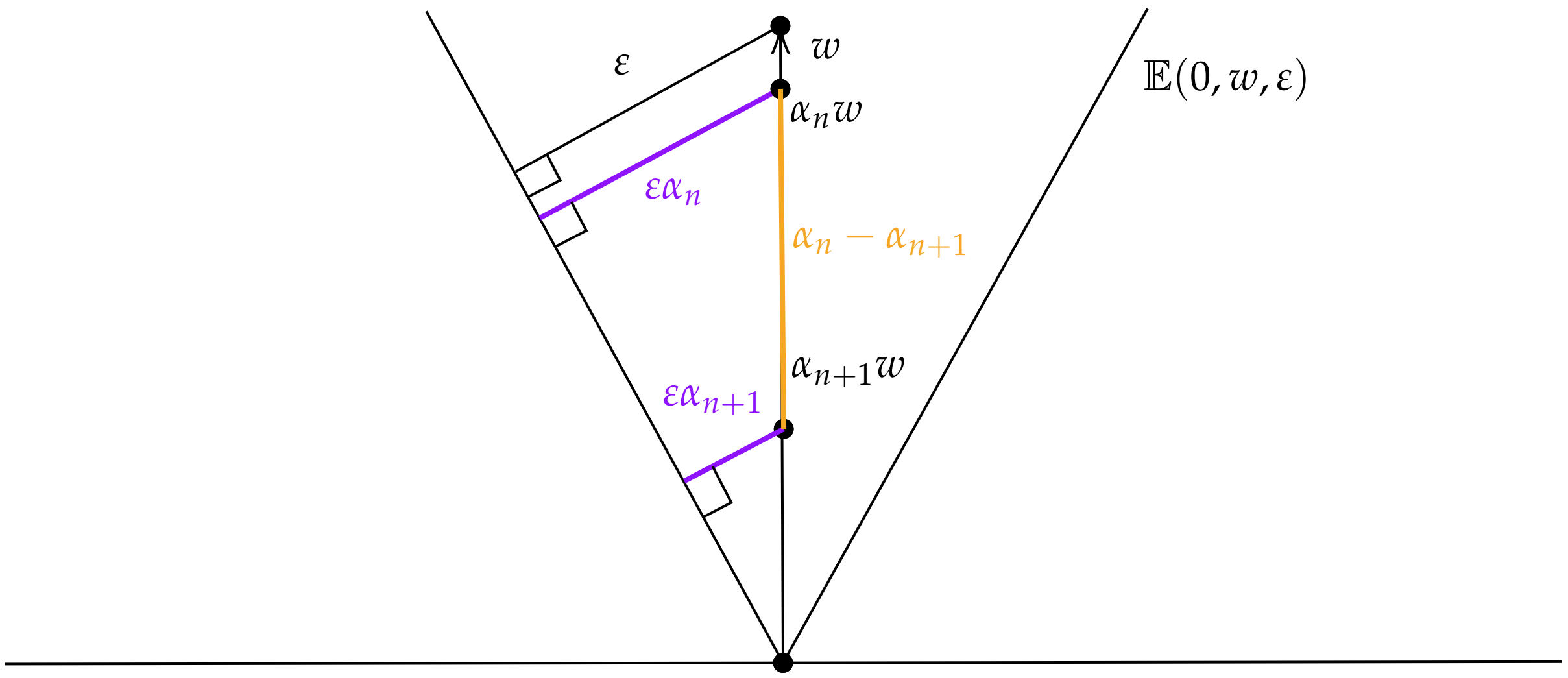

    \Needspace{6\baselineskip}
    \medskip\noindent\textbf{Situation 2. Points approaching direction $w$.}

    Again using the definition of a tangent vector, choose nonzero points $z_n\in g(B)$ such that
    \[
    \|z_n\|\longrightarrow0,
    \quad \frac{z_n}{\|z_n\|}\longrightarrow w.
    \]
    Indeed, the definition gives $\tau_n>0$ with $\tau_nz_n\to w$; then $\tau_n\|z_n\|\to\|w\|=1$, which proves the convergence of the directions. Passing to a subsequence, we may assume that
    $\|z_{n+1}\|\le\|z_n\|/2$ for every $n$.
    For all sufficiently large $n$,
    \[
    z_n\in g(B)\cap
       \mathbb E\left(0,w,\frac{\varepsilon}{2}\right).
    \]
    Join $z_n$ to $z_{n+1}$ by an arc in $g(B)$. As in Situation 1, follow it until it first reaches the next point, the boundary of the outer cone, or $\partial B(0,2\|z_n\|)$.
    
    The angle between a side of the outer cone and the closer side of the inner cone is
    \[
    \beta\coloneqq\arcsin\varepsilon
                    -\arcsin(\varepsilon/2)>0.
    \]
    This angle is independent of $n$. Using the formula 
    \[\sin\left(\arcsin\varepsilon-\arcsin(\varepsilon/2)\right)=\varepsilon/2\left(\sqrt{4-\varepsilon^2}-\sqrt{1-\varepsilon^2}\right)\]
    we obtain that $0<\sin\beta<\varepsilon<1/2$. Since $z_n$ lies in the inner cone, its distance from the boundary of the outer cone is at least $\|z_n\|\sin\beta$ (the boundary configuration is illustrated in Figure~\ref{f5}). Its distance from the outer sphere is $\|z_n\|$, whereas
    \[
    \|z_n-z_{n+1}\|\ge\|z_n\|-\|z_{n+1}\|
       \ge\frac12\|z_n\|.
    \]
    Thus the retained connected part has measure at least
    \[
    \min\left\{\frac12\|z_n\|,\|z_n\|\sin\beta,\|z_n\|\right\}
       =\|z_n\|\sin\beta.
    \]
    We obtain
    \[
    \hau^1\left(g(B)\cap\overline{B(0,2\|z_n\|)}
                    \cap\mathbb E(0,w,\varepsilon)\right)
       \ge\|z_n\|\sin\beta.
    \]
    Combining this with~\eqref{naef} at $R_1=2\|z_n\|$, for all sufficiently large $n$, gives
    \[
    \hau^1\left(g(B)\cap\overline{B(0,2\|z_n\|)}\right)
       \ge4\|z_n\|+\|z_n\|\sin\beta.
    \]
    Therefore
    \[
    \theta_1^*(g(B),0)
       \ge\limsup_{n\to\infty}
          \frac{4\|z_n\|+\|z_n\|\sin\beta}{4\|z_n\|}
       =1+\frac{\sin\beta}{4}>1,
    \]
    again contradicting~\eqref{zaloz}. 

\vspace{0.1cm}
\begin{figure}[H]
\centering

\tikzset{every picture/.style={line width=0.75pt}}

\begin{tikzpicture}[x=0.75pt,y=0.75pt,yscale=-1,xscale=1]

\draw    (331,252.23) -- (329.99,7.85) ;
\draw [shift={(329.99,5.85)}, rotate = 89.76] [color={rgb, 255:red, 0; green, 0; blue, 0 }  ][line width=0.75]    (10.93,-3.29) .. controls (6.95,-1.4) and (3.31,-0.3) .. (0,0) .. controls (3.31,0.3) and (6.95,1.4) .. (10.93,3.29)   ;
\draw    (227,-1.27) -- (331,252.23) ;
\draw [color={rgb, 255:red, 0; green, 0; blue, 0 }  ,draw opacity=1 ][line width=0.75]    (331,252.23) ;
\draw [shift={(331,252.23)}, rotate = 0] [color={rgb, 255:red, 0; green, 0; blue, 0 }  ,draw opacity=1 ][fill={rgb, 255:red, 0; green, 0; blue, 0 }  ,fill opacity=1 ][line width=0.75]      (0, 0) circle [x radius= 3.35, y radius= 3.35]   ;
\draw [shift={(331,252.23)}, rotate = 0] [color={rgb, 255:red, 0; green, 0; blue, 0 }  ,draw opacity=1 ][line width=0.75]      (0, 0) circle [x radius= 3.35, y radius= 3.35]   ;
\draw [color={rgb, 255:red, 0; green, 0; blue, 0 }  ,draw opacity=1 ][line width=0.75]    (329.99,5.85) ;
\draw [shift={(329.99,5.85)}, rotate = 0] [color={rgb, 255:red, 0; green, 0; blue, 0 }  ,draw opacity=1 ][fill={rgb, 255:red, 0; green, 0; blue, 0 }  ,fill opacity=1 ][line width=0.75]      (0, 0) circle [x radius= 3.35, y radius= 3.35]   ;
\draw [shift={(329.99,5.85)}, rotate = 0] [color={rgb, 255:red, 0; green, 0; blue, 0 }  ,draw opacity=1 ][line width=0.75]      (0, 0) circle [x radius= 3.35, y radius= 3.35]   ;
\draw [color={rgb, 255:red, 0; green, 0; blue, 0 }  ,draw opacity=1 ][line width=0.75]    (247,45.73) ;
\draw [shift={(247,45.73)}, rotate = 0] [color={rgb, 255:red, 0; green, 0; blue, 0 }  ,draw opacity=1 ][fill={rgb, 255:red, 0; green, 0; blue, 0 }  ,fill opacity=1 ][line width=0.75]      (0, 0) circle [x radius= 3.35, y radius= 3.35]   ;
\draw [shift={(247,45.73)}, rotate = 0] [color={rgb, 255:red, 0; green, 0; blue, 0 }  ,draw opacity=1 ][line width=0.75]      (0, 0) circle [x radius= 3.35, y radius= 3.35]   ;
\draw [color={rgb, 255:red, 0; green, 0; blue, 0 }  ,draw opacity=1 ][line width=0.75]    (296.28,166.76) ;
\draw [shift={(296.28,166.76)}, rotate = 0] [color={rgb, 255:red, 0; green, 0; blue, 0 }  ,draw opacity=1 ][fill={rgb, 255:red, 0; green, 0; blue, 0 }  ,fill opacity=1 ][line width=0.75]      (0, 0) circle [x radius= 3.35, y radius= 3.35]   ;
\draw [shift={(296.28,166.76)}, rotate = 0] [color={rgb, 255:red, 0; green, 0; blue, 0 }  ,draw opacity=1 ][line width=0.75]      (0, 0) circle [x radius= 3.35, y radius= 3.35]   ;
\draw [color={rgb, 255:red, 144; green, 19; blue, 254 }  ,draw opacity=1 ][line width=1.5]    (272.01,188.23) -- (296.28,166.76) ;
\draw [color={rgb, 255:red, 144; green, 19; blue, 254 }  ,draw opacity=1 ][line width=1.5]    (189,98.73) -- (247,45.73) ;
\draw   (284.05,178.14) -- (293.58,189.52) -- (281.53,199.61) ;
\draw   (199.12,89.93) -- (207.37,99.41) -- (197.26,108.22) ;
\draw    (29.88,252.68) -- (632.13,251.78) ;
\draw    (434,-2.27) -- (331,252.23) ;
\draw [color={rgb, 255:red, 245; green, 166; blue, 35 }  ,draw opacity=1 ][line width=1.5]    (247,45.73) -- (296.28,166.76) ;
\draw    (100,0.73) -- (331,252.23) ;
\draw    (560,-1.27) -- (331,252.23) ;

\draw (339.99,1.25) node [anchor=north west][inner sep=0.75pt]    {$\mathnormal{w}$};
\draw (255,25.63) node [anchor=north west][inner sep=0.75pt]    {$z_{n}$};
\draw (293,131.88) node [anchor=north west][inner sep=0.75pt]    {$z_{n+1}$};
\draw (558,11.63) node [anchor=north west][inner sep=0.75pt]    {$\mathbb{E}( 0,w,\varepsilon)$};
\draw (257.5,112) node [anchor=center,rotate=-67.85][inner sep=0.75pt]    {$\textcolor[rgb]{0.96,0.65,0.14}{|z_{n} |-|z_{n+1} |}$};
\draw (159,45.4) node [anchor=north west][inner sep=0.75pt]    {$\mathnormal{\textcolor[rgb]{0.56,0.07,1}{|z}\textcolor[rgb]{0.56,0.07,1}{_{n}}}\textcolor[rgb]{0.56,0.07,1}{|}\textcolor[rgb]{0.56,0.07,1}{\sin}\textcolor[rgb]{0.56,0.07,1}{\beta }$};
\draw (188,180.4) node [anchor=north west][inner sep=0.75pt]    {$\mathnormal{\textcolor[rgb]{0.56,0.07,1}{|z}\textcolor[rgb]{0.56,0.07,1}{_{n+1}}}\textcolor[rgb]{0.56,0.07,1}{|}\textcolor[rgb]{0.56,0.07,1}{\sin}\textcolor[rgb]{0.56,0.07,1}{\beta }$};
\draw (436,1.13) node [anchor=north west][inner sep=0.75pt]    {$\mathbb{E}\left( 0,w,\frac{\varepsilon}{2}\right)$};
\draw (302,204.4) node [anchor=north west][inner sep=0.75pt]    {$\beta $};

\end{tikzpicture}

\selectlanguage{english}
%
%
%
%
%
%
\caption{Situation 2: the boundary configuration giving the lower bound $|z_n|\sin\beta$. The two points shown lie on the same ray.}
\label{f5}
\end{figure}
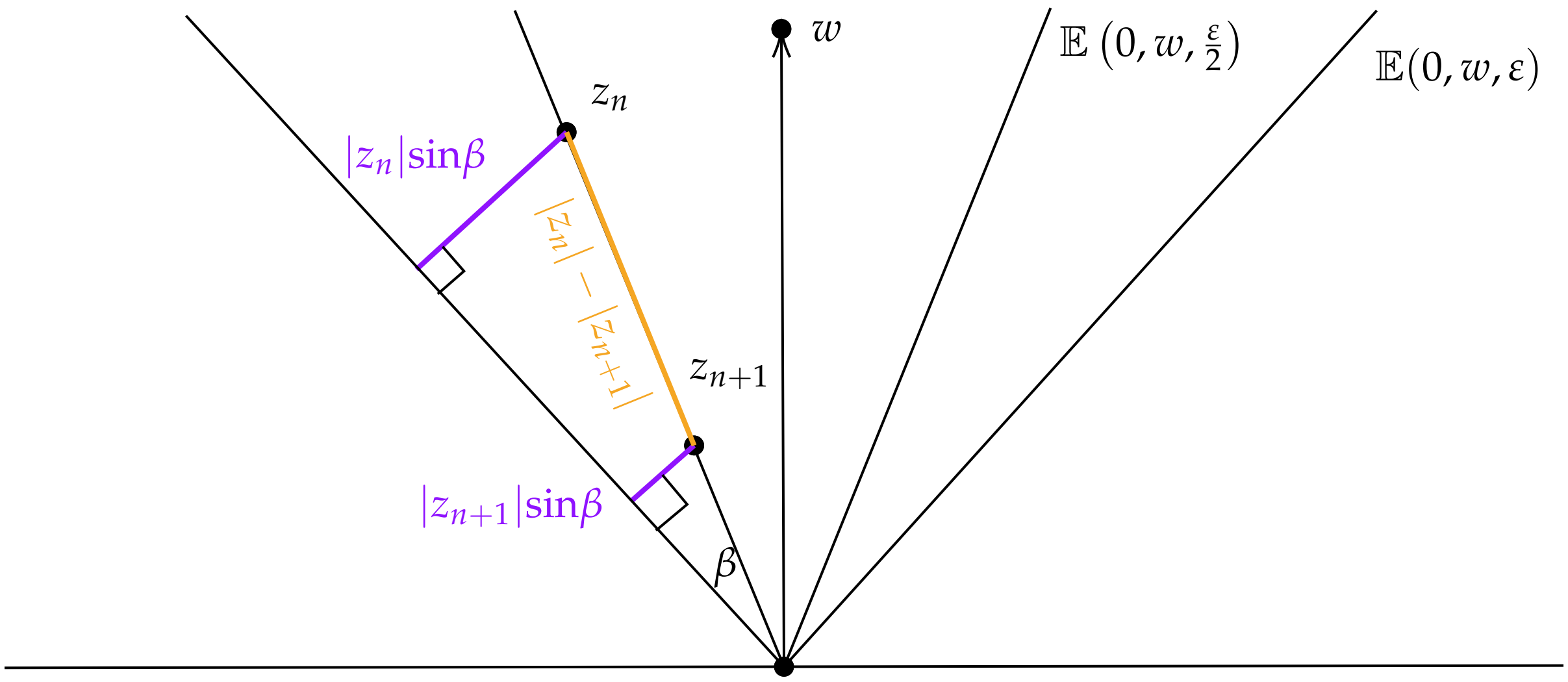
This completes the reverse inclusion and the proof of the theorem.
\end{proof.}




\section{An example without upper density bound and not satisfying the tangent-cone property}

Lipschitz continuity and the rank-one condition at a prescribed interior point of differentiability do not by themselves imply either the upper one-dimensional Hausdorff density bound by one at its image or the tangent-cone identity \eqref{imdg}. The following example demonstrates both failures, even though the differential has rank one at every interior point of differentiability.

Let $B=\overline{B(0,1)}\subset\R^2$ be the closed unit ball. Define $g\colon B\to\R^2$ by
\[
g(x,y)=
\begin{cases}
(x,0),
    & -1\le x\le \dfrac14;\\[6pt]
\left(\dfrac12-x,0\right),
    & \dfrac14<x\le \dfrac12;\\[6pt]
\left(0,x-\dfrac12\right),
    & \dfrac12<x\le 1.
\end{cases}
\]
We consider the point $(0,0)\in B(0,1)$.

\vspace{0.2cm}
Firstly, we verify that $g$ is Lipschitz. Evidently, $g$ is continuous on $B$. Let $-1\le s<t\le1$, and subdivide $[s,t]$ at whichever of the points $1/4$ and $1/2$ lie in its interior. Denote the resulting subdivision by
\[
s=t_0<t_1<\ldots<t_m=t.
\]
On each interval $[t_{j-1},t_j]$, the mapping $\tau\mapsto g(\tau,0)$ is affine, with slope equal to one of the unit vectors: $(1,0),\, (-1,0),\, (0,1)$. Hence, we can estimate
\[
\begin{aligned}
\|g(t,0)-g(s,0)\|
&\le
\sum_{j=1}^{m}
\|g(t_j,0)-g(t_{j-1},0)\|=
\sum_{j=1}^{m}(t_j-t_{j-1})
=t-s.
\end{aligned}
\]
The same inequality, with $|t-s|$ on the right, holds for all $s,t\in[-1,1]$ by symmetry. Since $g$ does not depend on its second coordinate, for all $(x,y),(x',y')\in B$ we obtain
\[
\begin{aligned}
\|g(x,y)-g(x',y')\|
&=\|g(x,0)-g(x',0)\|\le |x-x'|\le \|(x,y)-(x',y')\|.
\end{aligned}
\]
Therefore $g$ is $1$-Lipschitz.

\vspace{0.2cm}
Next, at every interior point with $x\notin\{1/4,1/2\}$, the mapping $g$ agrees with an affine mapping in a neighborhood of that point. Its differential is
\[
Dg(x,y)=
\begin{cases}
\begin{pmatrix}
1&0\\
0&0
\end{pmatrix},
    & x<\dfrac14,\\[12pt]
\begin{pmatrix}
-1&0\\
0&0
\end{pmatrix},
    & \dfrac14<x<\dfrac12,\\[12pt]
\begin{pmatrix}
0&0\\
1&0
\end{pmatrix},
    & x>\dfrac12.
\end{cases}
\]
All three matrices have rank one.

At an interior point with $x=1/4$, the left and right derivatives with respect to the first coordinate are $(1,0)$ and $(-1,0)$, respectively. At an interior point with $x=1/2$, they are $(-1,0)$ and $(0,1)$, respectively. Both one-sided variations of the first coordinate are admissible because the point is interior in $B$. Consequently, $g$ is not differentiable at any of these points.

Thus $Dg$ has rank one at every interior point of differentiability of $g$. In particular,
\[
g(0,0)=0,
\quad
Dg(0,0)=
\begin{pmatrix}
1&0\\
0&0
\end{pmatrix},
\quad
\operatorname{im}Dg(0,0)=\R\times\{0\}.
\]

\vspace{0.2cm}
We now describe the image of $g$. Every number in $[-1,1]$ is the first coordinate of a point of $B$, and $g$ depends only on this coordinate. The first branch of the definition covers the segment $[-1,1/4]\times\{0\}$. The second branch gives $[0,1/4]\times\{0\}$, while the third branch covers $\{0\}\times(0,1/2]$. Therefore
\[
E\coloneqq g(B)
=
\left([-1,\tfrac14]\times\{0\}\right)
\cup
\left(\{0\}\times[0,\tfrac12]\right).
\]

For every $0<r<1/4$, this gives the exact equality
\[
E\cap\overline{B(0,r)}
=
([-r,r]\times\{0\})
\cup
(\{0\}\times[0,r]).
\]
The two segments on the right have lengths $2r$ and $r$, respectively, and their intersection consists only of the origin. This yields
\[
\hau^1\left(E\cap\overline{B(0,r)}\right)
=2r+r=3r.
\]

We conclude that the upper density is given by
\[
\theta_1^*(g(B),g(0,0))
=
\limsup_{r\to0^+}
\frac{\hau^1\left(E\cap\overline{B(0,r)}\right)}{2r}
=
\frac32>1.
\]
Thus the upper density bound by $1$ fails at $g(0,0)=(0,0)$.

\vspace{0.2cm}
We next show that
\[
\operatorname{Tan}(E,0)
=
(\R\times\{0\})
\cup
(\{0\}\times[0,\infty)).
\]
We explicitly verify both inclusions using Definition~\ref{styczny}.

Let $v\in\operatorname{Tan}(E,0)$. For every positive integer $n$, there exist $z_n\in E$ and $r_n>0$ such that
\[
\|z_n\|<\frac1n,
\quad
\|r_nz_n-v\|<\frac1n.
\]
The closed set
\[
(\R\times\{0\})
\cup
(\{0\}\times[0,\infty))
\]
contains $E$ and is invariant under multiplication by positive scalars. Hence every $r_nz_n$ belongs to this set. Since $r_nz_n\to v$, we conclude that
\[
v\in
(\R\times\{0\})
\cup
(\{0\}\times[0,\infty)).
\]

Conversely, suppose that
\[
v\in
(\R\times\{0\})
\cup
(\{0\}\times[0,\infty)),
\]
and let $\varepsilon>0$. Choose $t>0$ such that
\[
t\|v\|<\min\{\varepsilon,1/4\}.
\]
By the description of $E$, the point $z=tv$ belongs to $E$: if $v$ is horizontal, then $z$ lies in $(-1/4,1/4)\times\{0\}$; if $v$ is vertical with nonnegative second coordinate, then $z$ lies in $\{0\}\times[0,1/4)$. This includes the case $v=0$.

Moreover,
\[
\|z\|<\varepsilon,
\quad
\left\|\frac1t z-v\right\|=0<\varepsilon.
\]
Taking $r=1/t>0$ in Definition~\ref{styczny}, we obtain $v\in\operatorname{Tan}(E,0)$. This proves the claimed description of the tangent cone.

In particular, $(0,1)\in\operatorname{Tan}(E,0)$, whereas $(0,1)\notin\operatorname{im}Dg(0,0)$. Consequently,
\[
\operatorname{im}Dg(0,0)
=
\R\times\{0\}
\subsetneq
\operatorname{Tan}(g(B),g(0,0)).
\]
We verified that the tangent-cone identity does not hold.

\vspace{0.2cm}
To sum up, we constructed a Lipschitz mapping whose differential has rank one at every interior point of differentiability, and which is differentiable at the interior point $(0,0)$, but with the following properties
\[
\begin{gathered}
\theta_1^*(g(B),g(0,0))=\frac32>1;\\[4pt]
\operatorname{im}Dg(0,0)
\subsetneq
\operatorname{Tan}(g(B),g(0,0)).
\end{gathered}
\]
Hence the assumptions preceding the density condition in Theorem~\ref{dodatek} imply neither the upper density bound nor the tangent-cone identity.

\vspace{0.2cm}
The example also satisfies $g(B)=g(\partial B)$. Indeed, for every $(x,y)\in B$, the point
\[
\left(x,\sqrt{1-x^2}\right)
\]
belongs to $\partial B$ and has the same image under $g$. This property was crucial in the rigidity proof in \cite{chlebik}.

\section*{Acknowledgements}
The author would like to thank Sławomir Kolasiński and Bernd Kirchheim for helpful discussions concerning the subject of this paper. The author is particularly grateful to Sławomir Kolasiński for pointing out the related result in Lemma~10.2 of Fang and Kolasiński~\cite{slawek}.

\section*{Funding}
This research was partially funded by the
\href{https://ncn.gov.pl/}{National Science Centre Poland}
grant number 2022/46/E/ST1/00328.

\section*{Disclaimer on the use of LLMs.}

During the preparation of this manuscript, the author used
GPT-6.1 Sol solely for language editing, stylistic revision, and for checking selected parts of the proofs for possible inconsistencies or gaps. The proof idea, scheme and all its underlying arguments were developed by the author, 
prior to the use of generative AI in the preparation of this manuscript. All mathematical arguments and figures are the author's own. The author takes entire responsibility for the final version and content of the manuscript.



\selectlanguage{english}


\begin{thebibliography}{99}

\bibitem{chlebik} M. Chlebík and B. Kirchheim, \textit{Rigidity for the four gradient problem}, Journal für die Reine und Angewandte Mathematik 551, 1--9 (2002).

%

\bibitem{slawek} Y. Fang and S. Kolasiński, \textit{Existence of solutions to a general geometric elliptic variational problem}. Calc. Var. 57, 91 (2018).
%
%

\bibitem{federer} H. Federer, \textit{Geometric Measure Theory}, Springer-Verlag (1969).


\end{thebibliography}
\end{document}